\documentstyle[11pt]{article}

\newtheorem{theo}{{\bf Theorem}}

\newtheorem{lemma}[theo]{{\bf Lemma}}
\newtheorem{prop}[theo]{{\bf Proposition}}
\newtheorem{coro}[theo]{{\bf Corollary}}
\newtheorem{defn}{{\bf Definition}}
\newtheorem{remark}{{\bf Remark}}

\font\bbb=msbm10 scaled\magstep1

\newcommand{\CC}{\mbox{\bbb C}}

\newcommand{\RR}{\mbox{\bbb R}}
\newcommand{\ZZ}{\mbox{\bbb Z}}

\def\C{{\kern.24em \vrule width.04em
  height1.1ex depth-.05ex \kern-.26em C}}

 \def\CCC{{\rm \kern.24em
 \vrule width.08em
   height1.5ex depth-.08ex
 \kern-.36em C}}

\newcommand{\circledast}{\odot \hspace{-3.33mm}\ast ~}

\newcommand{\si}{\sigma}

\begin{document}

\title{\bf From the icosahedron to natural triangulations of
\boldmath{$\CCC {\rm P}^{2}$} and \boldmath{${\rm S}^{2} \times
{\rm S}^{2}$}}
\author{\bf Bhaskar Bagchi {\boldmath $\cdot$} Basudeb Datta}

\date{}

\maketitle

%\vspace{-5mm}

\footnotetext[1]{The research of B. Datta was supported by
UGC grant UGC-SAP/DSA-IV.

\vspace{2mm}

Bhaskar Bagchi

Theoretical Statistics and Mathematics
Unit, Indian Statistical Institute,  Bangalore 560\,059, India.

e-mail: bbagchi@isibang.ac.in

\medskip

Basudeb Datta

Department of Mathematics, Indian Institute of
Science, Bangalore 560\,012,  India.

e-mail: dattab@math.iisc.ernet.in}

\begin{center}

\date{To appear in `{\bf Discrete \& Computational Geometry}'}

\end{center}

\smallskip

%\hrule

\bigskip

\noindent {\bf Abstract.} We present two constructions in this
paper\,: $(a)$ A 10-vertex triangulation $\CC P^{\,2}_{10}$ of the
complex projective plane $\CC P^{\,2}$ as a subcomplex of the join
of the standard sphere ($S^{\,2}_4$) and the standard real
projective plane ($\RR P^{\,2}_6$, the decahedron), its
automorphism group is $A_4$; $(b)$ a 12-vertex triangulation
$(S^{\,2} \times S^{\,2})_{12}$ of $S^{\,2} \times S^{\,2}$ with
automorphism group $2S_5$, the Schur double cover of the symmetric
group $S_5$. It is obtained by generalized bistellar moves from a
simplicial subdivision of the standard cell structure of $S^{\,2}
\times S^{\,2}$. Both constructions have surprising and intimate
relationships with the icosahedron. It is well known that
$\CC P^{\,2}$ has $S^{\,2}
\times S^{\,2}$ as a two-fold branched cover; we construct the
triangulation $\CC P^{\,2}_{10}$ of $\CC P^{\,2}$ by presenting a
simplicial realization of this covering map $S^{\,2} \times
S^{\,2} \to \CC P^{\,2}$. The domain of this simplicial map is a
simplicial subdivision of the standard cell structure of $S^{\,2}
\times S^{\,2}$, different from the triangulation alluded to in
$(b)$. This gives a new proof that K\"{u}hnel's $\CC P^{\,2}_9$
triangulates $\CC P^{\,2}$. It is also shown that $\CC
P^{\,2}_{10}$ and $(S^{\,2} \times S^{\,2})_{12}$ induce the
standard piecewise linear structure on $\CC P^{\,2}$ and $S^{\,2}
\times S^{\,2}$ respectively.

\bigskip

{\small


\smallskip


}

\bigskip

%\hrule

\section{Introduction and Results}

It is well known that the minimal triangulation $\RR P^{\,2}_6$ of
the real projective plane arises naturally from the icosahedron.
Indeed, it is the quotient of the boundary complex of the
icosahedron by the antipodal map. In this note, we report the
surprising result that there is a small triangulation (using only
10 vertices) of the complex projective plane which is also
intimately related to the icosahedron. Indeed, this simplicial
complex $\CC P^{\,2}_{10}$ occurs as a subcomplex of the
simplicial join $S^{\,2}_4\ast \RR P^{\,2}_6$. Our starting point
is the beautiful fact that $\CC P^{\,2}$ is homeomorphic to the
symmetrised square $(S^{\,2}\times S^{\,2})/\ZZ_2$ of the
2-sphere, where $\ZZ_2$ acts by co-ordinate flip. So, letting
$S^{\,2}_4$ denote the 4-vertex triangulation of $S^{\,2}$ (i.e.,
the boundary complex of the tetrahedron), we look for a
$\ZZ_2$-stable simplicial subdivision of the product cell complex
$S^{\,2}_4\times S^{\,2}_4$, without introducing extra vertices.
In order to ensure that the quotient complex (after quotienting by
$\ZZ_2$) does triangulate the quotient space $(S^{\,2}\times
S^{\,2})/\ZZ_2 = \CC P^2$, the $\ZZ_2$ action on this simplicial
subdivision must be ``pure" (cf. Definition \ref{pure} and Lemma
\ref{L1} in Section 5). It turns out that the following  $(S^{\,2}
\times S^{\,2})_{16}$ is the unique 16-vertex triangulation
satisfying these requirements.

\bigskip

\noindent {\bf Description of \boldmath{$(S^{\,2}\times
S^{\,2})_{16}$}\,:} The vertices are $x_{ij}$, $1\leq i, j \leq
4$. The full automorphism group is $A_4 \times \ZZ_2$, where $A_4$
acts on the indices and $\ZZ_2$ acts by $x_{ij} \leftrightarrow
x_{ji}$. Modulo this group the facets (maximal simplices) are the
following\,:
$$
x_{11}x_{22}x_{33}x_{12}x_{13}, ~ x_{11}x_{22}x_{12}x_{14}x_{34},
~ x_{11}x_{22}x_{14}x_{24}x_{34}, ~
x_{11}x_{22}x_{21}x_{24}x_{31}, ~ x_{11}x_{22}x_{24}x_{31}x_{34}.
$$
The full list of facets of $(S^{\,2}\times S^{\,2})_{16}$ may be
obtained from these five basic facets by applying the group
$A_4\times \ZZ_2$. Under this group, the first three basic facets
form orbits of length 24 each, while each of the last two forms an
orbit of length 12, yielding a total of $3\times 24 + 2 \times 12
= 96$ facets. It may be verified that the face vector of
$(S^{\,2}\times S^{\,2})_{16}$ is $(16, 84, 216, 240, 96)$.

\medskip

\noindent {\bf Description of \boldmath{$\CC P^{\,2}_{10}$}\,:}
Quotienting the above $(S^{\,2}\times S^{\,2})_{16}$ by the group
$\ZZ_2$ generated by the automorphism $x_{ij} \leftrightarrow
x_{ji}$, we get the $\CC P^{\,2}_{10}$ mentioned above. Its
vertices are $x_{ij}$, $1\leq i \leq j \leq 4$. Let $\alpha$,
$\beta$ be the generators of the alternating group $A_4$ given by
$\alpha = (1 2 3)$, $\beta = (12)(34)$. Then $\alpha$, $\beta$ act
on the vertices of $\CC P^{\,2}_{10}$ by\,:
\begin{eqnarray*}
\alpha \equiv
(x_{11}x_{22}x_{33})(x_{23}x_{13}x_{12})(x_{24}x_{34}x_{14}), &
\beta \equiv
(x_{11}x_{22})(x_{33}x_{44})(x_{24}x_{13})(x_{14}x_{23}).
\end{eqnarray*}
The following are the basic facets of $\CC P^{\,2}_{10}$ modulo
$A_4 = \langle\alpha, \beta\rangle$\,:
$$
x_{11}x_{22}x_{33}x_{12}x_{13}, ~ x_{11}x_{22}x_{12}x_{14}x_{34},
~ x_{11}x_{22}x_{14}x_{24}x_{34}, ~
x_{11}x_{22}x_{12}x_{13}x_{24}, ~ x_{11}x_{22}x_{13}x_{24}x_{34}.
$$
The full list of facets of $\CC P^{\,2}_{10}$ may be obtained from
these five basic facets by applying the group $A_4$. Under this
group, the first three basic facets form orbits of length 12 each,
while each of the last two forms an orbit of length 6, yielding a
total of $3\times 12 + 2 \times 6 = 48$ facets. The complex $\CC
P^{\,2}_{10}$ is 2-neighbourly and its face vector is $(10, 45,
110, 120, 48)$.

\bigskip

Here we prove the following\,:

\begin{theo}$\!\!\!${\bf .} \label{T1}
There are exactly two $16$-vertex simplicial complexes which $(i)$
are simplicial subdivisions of the cell complex $S^{\,2}_4
\times S^{\,2}_4$, $(ii)$ retain the self-homeomorphism $\alpha :
(x, y) \mapsto (y, x)$ of $|S^{\,2}_4| \times |S^{\,2}_4|$ as a
simplicial automorphism, and $(iii)$ the action of $\ZZ_2 =
\langle\alpha\rangle$ is pure $($cf. Definition $\ref{pure})$.
These two complexes are isomorphic and one of them is $(S^{\,2}
\times S^{\,2})_{16}$.
\end{theo}

\begin{coro}$\!\!\!${\bf .} \label{C2}
The complex $\CC P^{\,2}_{10} := (S^{\,2}\times
S^{\,2})_{16}/\ZZ_2$ is a $10$-vertex triangulation of $\CC
P^{\,2}$. Its full automorphism group is $A_4$.
\end{coro}

Let $T$ and ${\cal I}$ denote the solid tetrahedron and the
icosahedron in $\RR^{3}$ respectively. Thus, the cell complex
$S^{\,2}_4 \times S^{\,2}_4$ alluded to above is a subcomplex of
the boundary complex of the product polytope $T \times T$ in
$\RR^{\,6}$. Although we do not present the details in this paper,
Theorem \ref{T1} can be strengthened (following the same line of
arguments) to show that there is a unique simplicial subdivision
$S^{\,5}_{16}$ of the cell complex $\partial(T \times T)$ which is
$\ZZ_2$-stable with a pure $\ZZ_2$-action. To our utter
surprise, it turns out that as an abstract simplicial complex,
$S^{\,5}_{16}$ is isomorphic to the combinatorial join
$S^{\,2}_4\ast S^{\,2}_{12}$ of the boundary complexes of $T$ and
${\cal I}$ respectively.

\begin{remark}$\!\!\!${\bf .} \label{R1}
{\rm
This last fact has the following geometric interpretation. Let $T
\circledast {\cal I}$ denote the convex hull of $T \cup {\cal I}$,
where $T$ and ${\cal I}$ sit in two (three-dimensional) affine
subspaces of $\RR^{6}$ meeting at a point which is in the interior
of both polyhedra. Then $T \circledast {\cal I}$ is a simplicial
6-polytope and the boundary complex of this polytope is
combinatorially isomorphic to a simplicial subdivision of the
boundary complex of $T \times T$. This geometric result cries out
for a geometric explanation; but we have none.}
\end{remark}

By the construction, $(S^{\,2} \times S^{\,2})_{16}$ is a
subcomplex of  $S^{\,2}_4\ast S^{\,2}_{12}$. Since the decahedron
$\RR P^{\,2}_6$ is the quotient of $S^{\,2}_{12} = \partial {\cal
I}$ by $\ZZ_2$, and $\ZZ_2$ acts trivially on $S^{\,2}_4$ (the
latter being the combinatorial child of the ``diagonal" $S^{\,2}$
in $S^{\,2} \times S^{\,2}$, i.e., the $S^{\,2}_4$ in Figure 1),
on passing to the quotient, we find the surprising inclusion
$$
\CC P^{\,2}_{10} \subseteq S^{\,2}_4\ast \RR P^{\,2}_{6}.
$$
Indeed, $S^{\,2}_4$ and $\RR P^{\,2}_{6}$ occur as induced
subcomplexes of $\CC P^{\,2}_{10}$ on a complementary pair
of vertex sets.
Since both $S^{\,2}_4$ and $\RR P^{\,2}_{6}$ are classical
objects, and the combinatorial join is such a well known operation
on simplicial complexes, this inclusion says that $\CC
P^{\,2}_{10}$ was all along sitting there right before our eyes\,!

The number 10 obtained here is not optimal. It is well known (cf.
\cite{bk, am, bd1, bd3}) that any triangulation of $\CC P^{\,2}$
requires at least nine vertices, and there is a unique 9-vertex
triangulation $\CC P^{\,2}_9$ of this manifold, obtained by
K\"{u}hnel (\cite{kb2, kl}). But, our construction is natural in
that it is obtained by a combinatorial mimicry of a topological
construction of $\CC P^{\,2}$. It shares this naturalness with
another 10-vertex triangulation, say $K^4_{10}$, of $\CC P^{\,2}$
available in the literature, namely the ``equilibrium"
triangulation of Banchoff and K\"{u}hnel (\cite{kb1}). Here we
prove the following\,:

\begin{theo}$\!\!\!${\bf .} \label{T3}
The simplicial complex $\CC P^{\,2}_{10}$ is bistellar equivalent
to both $\CC P^{\,2}_9$ and $K^4_{10}$.
\end{theo}

\begin{coro}$\!\!\!${\bf .} \label{C4}
$(a)$ K\"{u}hnel's $9$-vertex simplicial complex $\CC P^{\,2}_{9}$
triangulates $\CC P^{2}$. $(b)$ Both $\CC P^{\,2}_{10}$ and
$\CC P^{\,2}_{9}$ induce the standard pl-structure on $\CC P^{\,2}$.
\end{coro}

Of course, in principle these ideas generalize to arbitrary
dimensions. In general, the $d$-dimensional complex projective
space $\CC P^{\,d}$ is the symmetric $d$-th power of $S^{\,2}$,
i.e., the quotient of $(S^{\,2})^d$ by the symmetric group $S_d$
acting by co-ordinate permutations. Unfortunately, even in the
next case $d=3$, it is not possible to subdivide the cell complex
$S^{\,2}_4 \times S^{\,2}_4 \times S^{\,2}_4$ into a simplicial
complex, with a pure $S_3$-action, without adding more vertices.
Indeed, we found that we need to add 60 more vertices to obtain an
$(S^{\,2} \times S^{\,2} \times S^{\,2})_{124}$. On quotienting,
we obtain a $\CC P^{\,3}_{30}$ - again with full automorphism
group $A_4$. The details are so complicated that we decided to
postpone publication. We are presently trying to see if one can
apply bistellar moves to this $\CC P^{\,3}_{30}$ to reduce the
number of vertices. It is known that any triangulation of $\CC
P^{\,3}$  requires at least 17 vertices (cf. \cite{am}).

After we submitted a preliminary version of this paper to arXiv
(arXiv:1004.3157v1, 2010), Ulrich Brehm (\cite{br}) communicated
to us that he had
the idea of obtaining $\CC P^{\,2}_{10}$ as a quotient of a
16-vertex $S^{\,2} \times S^{\,2}$ in the 1980's; however he
never published the details.

\medskip

We obtain a second simplicial subdivision $(S^{\,2} \times
S^{\,2})_{16}^{\prime}$ of $S^{\,2}_4 \times S^{\,2}_4$.

\bigskip

\noindent {\bf Description of {\boldmath $(S^{\,2} \times
S^{\,2})_{16}^{\prime}$}\,:} This is a second simplicial
subdivision of the cell complex $S^{\,2}_4 \times S^{\,2}_4$. It
has the same vertex-set and automorphism group $A_4$. Modulo the
group $A_4$, its basic facets are\,:
\begin{eqnarray*}
& x_{11}x_{12}x_{13}x_{21}x_{31}, ~
x_{11}x_{12}x_{14}x_{21}x_{31}, ~ x_{11}x_{13}x_{14}x_{21}x_{31},
~ x_{12}x_{13}x_{23}x_{31}x_{32}, & \\
& x_{12}x_{14}x_{21}x_{24}x_{31}, ~
x_{12}x_{14}x_{24}x_{31}x_{34}, ~ x_{12}x_{21}x_{24}x_{31}x_{32},
~ x_{12}x_{24}x_{31}x_{32}x_{34}. &
\end{eqnarray*}
Each facets is in an orbit of length 12, yielding a total of $8
\times 12 = 96$ facets. The complex $(S^{\,2} \times
S^{\,2})_{16}^{\prime}$ has the same face vector as $(S^{\,2}
\times S^{\,2})_{16}$, namely, $(16, 84, 216, 240, 96)$.

\bigskip

We perform a finite sequence of generalized bistellar moves on
$(S^{\,2} \times S^{\,2})_{16}^{\prime}$ and obtain the
following 12-vertex triangulation $(S^{\,2} \times S^{\,2})_{12}$
of $S^{\,2} \times S^{\,2}$.

\bigskip

\noindent {\bf Description of {\boldmath $(S^{\,2} \times
S^{\,2})_{12}$}\,:} The vertices are $x_{ij}$, $1\leq i \neq j
\leq 4$. Its automorphism group $2S_5$ is generated by the two
automorphisms $h = (x_{12} x_{14} x_{21} x_{24} x_{31})
(x_{13} x_{42} x_{43} x_{32} x_{34})$ and $g = (x_{12} x_{21} x_{24}
x_{42} x_{14} x_{41} x_{43} x_{34} x_{13} x_{31} x_{32} x_{23})$.
Modulo this group, $(S^{\,2} \times
S^{\,2})_{12}$ is generated by the following two basic facets:
$$x_{12} x_{14} x_{21} x_{24}x_{31}, x_{12}x_{13}x_{14}x_{21}
x_{31}.
$$
The first basic facet is in an orbit of size 12, while
the second is in an orbit of size 60, yielding a total of 72
facets. Its face vector is $(12, 60, 160, 180, 72)$.

\begin{theo}$\!\!\!${\bf .} \label{T5}
The simplicial complex $(S^{\,2}\times S^{\,2})_{12}$ is a
triangulation of $S^{\,2}\times S^{\,2}$. Its full automorphism
group is $2S_5$, the non-split extension of $\ZZ_2$ by $S_5$.
\end{theo}

The complex $(S^{\,2} \times S^{\,2})_{12}$  has many remarkable
properties. Its automorphism group is transitive on
its vertices and edges. All its vertices have degree 10 and all
its edges have degree 8. Indeed, the link of each edge is
isomorphic to the 2-sphere $S^{\,2}_8$ obtained from the boundary
complex of the octahedron by starring two vertices in a pair of
opposite faces. Also, all triangles of $(S^{\,2} \times
S^{\,2})_{12}$ are of degree 3 or 5. The automorphism group is
transitive on its triangles of each degree. The degree 3 triangles
constitute a weak pseudomanifold whose strong components are two
icosahedra. Thus, we find a pair $I_1$, $I_2$ of icosahedra
sitting canonically inside the 2-skeleton of $(S^{\,2} \times
S^{\,2})_{12}$. These two icosahedra are ``antimorphic" in the
sense that the identity map is an antimorphism between them (cf.
Definition \ref{D1} below). The structure of $(S^{\,2} \times
S^{\,2})_{12}$ is completely described in terms of this
antimorphic pair of icosahedra. The full automorphism group $2S_5$
of $(S^{\,2} \times S^{\,2})_{12}$ is a double cover of the common
automorphism group of these two icosahedra.

Again, the number 12 here is not optimal. In \cite{kl}, K\"{u}hnel
and La{\ss}mann have shown that any triangulation of $S^{\,2}
\times S^{\,2}$ needs at least 11 vertices, and in \cite{lu}, Lutz
finds (via computer search) several 11-vertex triangulations of
$S^{\,2} \times S^{\,2}$, all with trivial automorphism groups.
Surprisingly, even though $(S^{\,2} \times S^{\,2})_{12}$ is not
minimal, it does not admit any proper bistellar moves. Thus, there
is no straightforward way to obtain a minimal triangulation of
$S^{\,2} \times S^{\,2}$ starting from $(S^{\,2} \times
S^{\,2})_{12}$.

In \cite{sp}, Sparla proved two remarkable inequalities on the
Euler characteristic $\chi$ of a combinatorial 4-manifold $M$
satisfying certain conditions. His first result is that if there
is a centrally symmetric simplicial polytope $P$ of dimension
$d\geq 6$ such that $M \subseteq \partial P$ and ${\rm skel}_2(M)
= {\rm skel}_2(P)$, then $10(\chi -2) \geq 4^3 {(d-1)/2 \choose
3}$. Equality holds here if and only if $P$ is a cross polytope
(i.e., dual of a hypercube). His second result is\,: if $M$ has
$2d$ vertices and admits a fixed point free involution then
$10(\chi -2) \leq 4^3 {(d-1)/2 \choose 3}$. Equality holds if and
only if $M$ also satisfies the hypothesis of the first result for
a cross polytope $P$. Notice that, in view of the Dehn-Sommerville
equations, equality in either inequality determines the face
vector of $M$ in terms of $d$ alone. To obtain an example of
equality (in both results) with $d=6$, Sparla searched for (and
found) a 4-manifold with the predicted face vector under the
assumption of an automorphism group $A_5 \times \ZZ_2$. To
determine the topological type of the resulting 12-vertex
4-manifold, he had to compute its intersection form and then
appeal to Freedman's classification of simply connected smooth
4-manifolds. We believe that our approach to Sparla's complex not
only elucidates its true genesis, but also reveals its rich
combinatorial structure and contributes to an elementary
determination of its topological type. Note, however, that
Sparla's approach reveals yet another remarkable property of
$(S^{\,2} \times S^{\,2})_{12}$. It provides a tight rectilinear
embedding of $S^{\,2} \times S^{\,2}$ in $\RR^6$.

\begin{remark}$\!\!\!${\bf .} \label{R2}
{\rm If $X$ is a triangulated $4$-manifold on at most $12$
vertices, then its vertex-links are homology $3$-spheres on at
most $11$ vertices, and hence (cf. \cite{bd5}) are combinatorial
spheres. Thus all triangulated $4$-manifolds on at most $12$
vertices are combinatorial manifolds. (More generally, this
argument yields\,: All triangulated $d$-manifolds on at most $d+8$
vertices are combinatorial manifolds.) In particular, both $\CC
P^{\,2}_{10}$ and $(S^{\,2}\times S^{\,2})_{12}$ are combinatorial
manifolds. Actually, an old result of Bing (\cite{bi}) says that
all the vertex links of any triangulated $4$-manifold are simply
connected triangulated $3$-manifolds. Therefore, in view of
Perelman's theorem (Poincar\'{e} conjecture) (\cite{pe}),  all
triangulated $4$-manifolds are combinatorial manifolds,
irrespective of the number of vertices.
}
\end{remark}

\begin{remark}$\!\!\!${\bf .} \label{R3}
{\rm In \cite{ap}, Akhmedov and Park have shown that $S^{\,2}
\times S^{\,2}$ has countably infinite number of distinct smooth
structures. Since there is an one to one correspondence between
the smooth structures and pl-structures on a 4-manifold (cf.
\cite[page 167]{sa}), it follows that $S^{\,2} \times S^{\,2}$ has
infinitely many distinct pl-structures. Since $(S^{\,2} \times
S^{\,2})_{16}$ and $(S^{\,2} \times S^{\,2})_{16}^{\prime}$ are
simplicial subdivisions of $S^{\,2}_4 \times S^{\,2}_4$, it
follows that the pl-structures given by $(S^{\,2} \times
S^{\,2})_{16}$ and $(S^{\,2} \times S^{\,2})_{16}^{\prime}$ are
standard. Again, $(S^{\,2} \times S^{\,2})_{12}$ is
combinatorially equivalent to $(S^{\,2} \times
S^{\,2})_{16}^{\prime}$ (cf. Remark \ref{R5}) and hence gives the
same pl-structure as $(S^{\,2} \times S^{\,2})_{16}^{\prime}$. So,
all the triangulations of $S^{\,2}\times S^{\,2}$ discussed here
give the standard pl-structure on $S^{\,2}\times S^{\,2}$. }
\end{remark}

\section{Preliminaries}

All simplicial complexes considered here are finite and the empty
set is a simplex (of dimension $-1$) of every simplicial complex.
We now recall some definitions here.

For a finite set $V$ with $d+2$ ($d\geq 0$) elements, the set
$\partial V$ (respectively, $\bar{V}$) of all the proper (resp.
all the) subsets of $V$ is a simplicial complex and triangulates
the $d$-sphere $S^{\,d}$ (resp. the $(d+1)$-ball). The complex
$\partial V$ is called the {\em standard $d$-sphere} and is also
denoted by $S^{\,d}_{d+2}(V)$ (or simply by $S^{\,d}_{d+2}$). The
complex $\bar{V}$ is called the {\em standard} $(d+1)$-ball and is
also denoted by $D^{d+1}_{d+2}(V)$ (or simply by $D^{d+1}_{d+2}$).
(Generally, we write $X = X^d_n$
to indicate that $X$ has $n$ vertices and dimension $d$.)

For simplicial complexes $X$, $Y$ with disjoint vertex-sets, their
{\em join} $X \ast Y$ is the simplicial complex whose simplices
are all the disjoint unions $A \cup B$ with $A\in X$, $B\in Y$.

If $\si$ is a simplex of a simplicial complex $X$ then the {\em
link} of $\si$ in $X$, denoted by ${\rm lk}_X(\sigma)$, is the
simplicial complex whose simplices are the simplices $\tau$ of $X$
such that $\tau \cap \si = \emptyset$ and $\si\cup\tau$ is a
simplex of $X$. The number of vertices in the link of $\si$ is
called the {\em degree} of $\si$. Also, the {\em star} of
$\sigma$, denoted by ${\rm star}_X(\sigma)$ or ${\rm
star}(\sigma)$, is the subcomplex $\bar{\sigma} \ast {\rm
lk}_X(\sigma)$ of $X$.

For a simplicial complex $X$, $|X|$ denotes the {\em geometric
carrier}. It may be described as the subspace of $[0, 1]^{V(X)}$
(where $V(X)$ is the vertex set of $X$) consisting of all
functions $f \colon V(X) \to [0, 1]$ satisfying (i) ${\rm
Support}(f) \in X$ and (ii) $\sum_{x\in V(X)} f(x) = 1$. If a
space $Y$ is homeomorphic to $|X|$ then we say that $X$ {\em
triangulates} $Y$. If $|X|$ is a topological manifold
(respectively, $d$-sphere) then $X$ is called a {\em triangulated
manifold} (resp. {\em triangulated $d$-sphere}). If $|X|$ is a pl
manifold (with the pl structure induced by $X$) then $X$ is called
a {\em combinatorial manifold}. For $1\leq d\leq 4$, $X$ is a
combinatorial $d$-manifold if and only if the vertex links are
triangulated $(d-1)$-spheres.

The {\em face vector} of a $d$-dimensional simplicial complex is
the vector $(f_0, f_1, \dots, f_d)$, where $f_i$ is the number of
$i$-dimensional simplices in the complex.

If $X$ is a $d$-dimensional pure simplicial complex (i.e., every
maximal simplex is $d$-dimensional) and $D$, $\widehat{D}$ are
triangulations of the $d$-ball such that (i) $\partial D =
\partial \widehat{D} =\widehat{D} \cap X$, and (ii) $D
\subseteq X$, then the simplicial complex $\widehat{X} := (X
\setminus D) \cup \tilde{D}$ is said to be obtained from $X$ by a
{\em generalized bistellar move} (GBM) with respect to the pair
$(D,\widehat{D})$. Clearly, in this case, $\widehat{X}$ and $X$
triangulate the same topological space and if $u$ is a vertex in
$\partial D$ then ${\rm lk}_{\widehat{X}}(u)$ is obtained from
${\rm lk}_{X}(u)$ by a GBM (cf. \cite{bd9}).

In particular, let $A$ be a simplex of $X$ whose link in $X$ is a
standard sphere $\partial B$. Suppose also that $B \not\in X$.
Then, we may perform the GBM with respect to the pair of balls
$(A\ast\partial B, B \ast \partial A)$. Such an operation is
called a {\em bistellar move}, and will be denoted by $A \mapsto
B$. Also, if $C$ is any simplex of $X$ and $x$ is a new symbol,
then we may perform the GBM on $X$ with respect to the pair
$(\bar{C} \ast {\rm lk}_X(C), (\bar{x} \ast \partial C) \ast {\rm
lk}_X(C))$. The resulting simplicial complex $\widehat{X}$ is
said to be obtained from $X$ by {\em starring} the vertex $x$ in
the simplex $C$. In case $C$ is a facet, this is a bistellar move
- the only sort of bistellar move which increases the number of
vertices. All other kinds of bistellar moves are said to be {\em
proper}.

Two pure simplicial complexes are called {\em bistellar equivalent}
if one is obtained from the other by a finite sequence of
bistellar moves. If $X$ is obtained from $Y$ by the bistellar
move $A \mapsto B$ then the complex $Z$ obtained from $Y$ by starring
a new vertex $u$ in $B$ is a subdivision of both $X$ and $Y$.
This implies that bistellar equivalent complexes induce same
pl-structure on their common geometric carrier.

The group $\ZZ_2$ acts on $S^{\,2} \times S^{\,2}$ by co-ordinate
flip. The following proposition is well known to algebraic
geometers (cf. \cite{la})\,:

\begin{prop}$\!\!\!${\bf .} \label{P6}
The quotient space $(S^{\,2} \times S^{\,2})/\ZZ_2$ is homeomorphic
to the complex projective plane $\CC P^{\,2}$.
\end{prop}

\section{Relations with the icosahedron}

\noindent {\bf Emergence of the icosahedron\,:} Let $T_0$ be the
tetrahedron with vertex-set $V = \{x_1, x_2$, $x_3, x_4\}$. Then,
viewed abstractly, the boundary complex of the product polytope
$T_0 \times T_0$ has vertex-set $V\times V$, and faces $A\times
B$, where $A$ and $B$ range over all the subsets of $V$. The
product cell complex for $S^{\,2}_4 \times S^{\,2}_4 = (\partial
T_0) \times (\partial T_0)$ is the subcomplex consisting of cells
$A\times B$, where $A$ and $B$ range over all the proper subsets
of $V$. We use the notation $x_{ij}$ to denote the vertex $(x_i,
x_j)$ of $T_0 \times T_0$. For $i \neq j$, $k \neq l$, $x_{ij}
x_{kl}$ forms an edge of $T_0 \times T_0$ if and only if it is one
of the solid edges of the icosahedron in Figure 1. (This picture
is a Schlegel diagram obtained by projecting the boundary of the
icosahedron on one of its faces. Thus, there is only one ``hidden"
face (namely, $x_{41}x_{42}x_{43}$) in the picture. What is
important for us is the label given to the vertices.)

Notice that the broken edges in the icosahedron are precisely the
edges $x_{ij}x_{kl}$ where $\{i, j, k, l\}$ is an even permutation
of $\{1, 2, 3, 4\}$.

To obtain the appropriate triangulation of $S^{\,2} \times
S^{\,2}$, we join $x_{ii}$ to all vertices for all $i$ and also
introduce the broken edges of the icosahedron. Thus viewed, one
sees the simplicial subdivision $(S^{\,2} \times S^{\,2})_{16}$ of
the cell complex $(\partial T_0) \times (\partial T_0)$ as a
subcomplex of $(\partial T) \ast  (\partial {\cal I})$, where $T$
is the tetrahedron with vertex-set $\{x_{ii} \, : \, 1\leq i\leq
4\}$ and ${\cal I}$ is the icosahedron depicted in Figure 1.

Notice also that the $\ZZ_2$-action $x_{ij}\leftrightarrow x_{ji}$
fixes the vertices of $T$ and acts on ${\cal I}$ as the antipodal
map. Thus, going modulo $\ZZ_2$, we find $\CC P^{\,2}_{10}$ as a
subcomplex of the 5-dimensional simplicial complex $S^{\,2}_4 \ast
\RR P^{\,2}_6$, where $S^{\,2}_4$ is the 4-vertex 2-sphere given
by the boundary complex of $T$ and $\RR P^{\,2}_6$ is the
(minimal) triangulation of the real projective plane (with
vertices of the same name being identified) given in Figure 1.

\medskip

%\hrule

\setlength{\unitlength}{1.75mm}

\begin{picture}(83,49.5)(-8,-6.5)

%%%%%%%%%%%%%%%%%%%%%%% %%%%%%%%%%%%%%%%%%%%%%%%%%%%%
{\boldmath{

\thicklines

\put(0,0){\line(1,0){60}} \put(0,0){\line(3,4){30}}
\put(0,0){\line(5,1){30}} \put(0,0){\line(5,3){20}}
\put(60,0){\line(-3,4){30}} \put(60,0){\line(-5,3){20}}
\put(60,0){\line(-1,1){20}} \put(20,20){\line(1,2){10}}
\put(30,24){\line(0,1){16}}

\put(20,12){\line(5,-3){10}} \put(20,12){\line(5,6){5}}
\put(20,12){\line(0,1){8}}

\put(30,12){\line(1,0){10}} \put(30,12){\line(-5,6){5}}
\put(30,12){\line(5,6){5}}

\put(30,6){\line(5,3){10}} \put(30,6){\line(0,1){6}}

\put(40,12){\line(0,1){8}}

\put(20,20){\line(5,-2){5}} \put(20,20){\line(5,2){10}}
\put(40,20){\line(-5,-2){5}} \put(40,20){\line(-5,2){10}}

\put(25,18){\line(1,0){10}} \put(35,18){\line(-5,6){5}}

}}

\thinlines

\put(-10.5,-6.5){\line(1,0){85}} \put(-10.5,43){\line(1,0){85}}
%%These two are boundary lines

%\put(0,0){\line(1,1){20}}
\put(0,0){\line(1,1){2.1}} \put(2.55,2.55){\line(1,1){2.1}}
\put(5.1,5.1){\line(1,1){2.1}} \put(7.65,7.65){\line(1,1){2.1}}
\put(10.2,10.2){\line(1,1){2.1}}
\put(12.75,12.75){\line(1,1){2.1}}
\put(15.3,15.3){\line(1,1){2.1}}
\put(17.85,17.95){\line(1,1){2.15}}

%\put(60,0){\line(-5,1){30}}
\put(30,6){\line(5,-1){3}} \put(33.5,5.3){\line(5,-1){3}}
\put(37,4.6){\line(5,-1){3}} \put(40.5,3.9){\line(5,-1){3}}
\put(44,3.2){\line(5,-1){3}} \put(47.5,2.5){\line(5,-1){3}}
\put(51,1.8){\line(5,-1){3}} \put(54.5,1.1){\line(5,-1){3}}
\put(60,0){\line(-5,1){2.1}}

%\put(40,20){\line(-1,2){10}}
\put(40,20){\line(-1,2){1.5}} \put(38.3,23.4){\line(-1,2){1.5}}
\put(36.6,26.8){\line(-1,2){1.5}}
\put(34.9,30.2){\line(-1,2){1.5}}
\put(33.2,33.6){\line(-1,2){1.5}} \put(31.5,37){\line(-1,2){1.5}}
%%%

%\put(20,12){\line(1,0){10}}
\put(20,12){\line(1,0){2.25}} \put(22.75,12){\line(1,0){2}}
\put(25.25,12){\line(1,0){2}} \put(27.75,12){\line(1,0){2.25}}

%\put(40,12){\line(-5,6){5}}
\put(40,12){\line(-5,6){1.68}} \put(38,14.4){\line(-5,6){1.68}}
\put(36,16.8){\line(-5,6){1.68}}

%\put(25,18){\line(5,6){5}}
\put(30,24){\line(-5,-6){1.68}} \put(28,21.6){\line(-5,-6){1.68}}
\put(26,19.2){\line(-5,-6){1.68}}

\put(-3,1){\mbox{$x_{43}$}} \put(60.5,1){\mbox{$x_{42}$}}
\put(32,39.5){\mbox{$x_{41}$}} \put(31,24.5){\mbox{$x_{31}$}}
\put(31,10){\mbox{$x_{14}$}} \put(29,4){\mbox{$x_{13}$}}
\put(41,12.5){\mbox{$x_{12}$}} \put(41,20){\mbox{$x_{32}$}}
\put(16.5,12.5){\mbox{$x_{23}$}} \put(16.5,20){\mbox{$x_{21}$}}
\put(27,18.75){\mbox{$x_{24}$}} \put(34.4,19.7){\mbox{$x_{34}$}}

\put(5,-4){\mbox{$\partial {\cal I} = S^{\,2}_{12}$}}

\put(30,-5){\mbox{\bf Figure 1}}

\thicklines

\put(50,28){\line(5,-3){10}} \put(50,28){\line(5,6){5}}
\put(50,28){\line(0,1){8}}

\put(60,28){\line(1,0){10}} \put(60,28){\line(-5,6){5}}
\put(60,28){\line(5,6){5}}

\put(60,22){\line(5,3){10}} \put(60,22){\line(0,1){6}}

\put(70,28){\line(0,1){8}}

\put(50,36){\line(5,-2){5}} \put(50,36){\line(5,2){10}}
\put(70,36){\line(-5,-2){5}} \put(70,36){\line(-5,2){10}}

\put(55,34){\line(1,0){10}} \put(65,34){\line(-5,6){5}}

%%%

\put(50,28){\line(1,0){10}} \put(70,28){\line(-5,6){5}}
\put(55,34){\line(5,6){5}}

\put(61,40.5){\mbox{$x_{13}$}} \put(61,26){\mbox{$x_{14}$}}
\put(57,20.2){\mbox{$x_{13}$}} \put(71,28.5){\mbox{$x_{12}$}}
\put(71,36){\mbox{$x_{23}$}} \put(46.5,28.5){\mbox{$x_{23}$}}
\put(46.5,36){\mbox{$x_{12}$}} \put(57,34.75){\mbox{$x_{24}$}}
\put(64.4,35.7){\mbox{$x_{34}$}}

\put(62,17){\mbox{$\RR P^{\,2}_{6}$}}

\thicklines

\put(-6,38){\line(1,0){20}} \put(4,23){\line(2,3){10}}
\put(4,23){\line(-2,3){10}} \put(4,23){\line(0,1){9}}
\put(4,32){\line(5,3){10}} \put(4,32){\line(-5,3){10}}

\put(5,22){\mbox{$x_{11}$}} \put(5,31){\mbox{$x_{44}$}}
\put(-7,39.4){\mbox{$x_{22}$}} \put(13,39.4){\mbox{$x_{33}$}}

\put(-4,19){\mbox{$\partial T = S^{\,2}_{4}$}}

\end{picture}

%\hrule

\medskip

From our nomenclature for the vertices, the inclusion $\CC
P^{\,2}_{10}\subseteq S^{\,2}_4 \ast \RR P^{\,2}_6$ is obvious, as
is the fact that $(\partial T) \ast (\partial {\cal I})$ is a
simplicial subdivision of the boundary complex of $T_0\times T_0$.

Finally, note that $\Delta_i = \{x_{ij} \, : \, j\neq i\}$ and
$\Delta^i = \{x_{ji} \, : \, j\neq i\}$ are triangles of the
icosahedron, and $\{\Delta_1,  \Delta_2,  \Delta_3,  \Delta_4\}$
and $\{\Delta^1, \Delta^2, \Delta^3, \Delta^4\}$ are antipodal
pairs of quadruples (consisting of triangles) partitioning the
vertex-set of the icosahedron. It is easy to see that there are
exactly five such pairs in the icosahedron, and the automorphism
group $A_5 \times \ZZ_2$ of ${\cal I}$ acts transitively on them.
The stabilizer of each such pair is $A_4\times \ZZ_2$, and $A_4$
acts regularly on the vertex-set of ${\cal I}$. Our choice of
nomenclature for the vertices of ${\cal I}$ amounts to choosing
one such antipodal pair of quadruples. This is because we have
$\Delta_i \cap \Delta^j = \emptyset$ if $i=j$ and $=\{x_{ij}\}$ if
$i\neq j$. Viewed dually, one sees Kepler's regular tetrahedra
embedded in the dodecahedron. Namely, the centres of $\Delta_i$,
$1\leq i\leq 4$ (as well as of $\Delta^i$, $1\leq i\leq 4$) are
the vertices of a regular tetrahedron inscribed in the dual
dodecahedron.

\bigskip

The 12-vertex triangulation $(S^{\,2} \times S^{\,2})_{12}$ of
$S^{\,2} \times S^{\,2}$ is obtained from $(S^{\,2} \times
S^{\,2})_{16}^{\prime}$ by a sequence of bistellar moves (cf.
proof of Theorem \ref{T5}). However, its most elegant description
requires the introduction of the following definition.

\begin{defn}$\!\!\!${\bf .} \label{D1}
{\rm Let $I_1$ and $I_2$ be two copies of the icosahedron. A
bijection $f \colon V(I_1) \to V(I_2)$ is said to be an {\em
antimorphism} if, for all vertices $x$, $y$ of $I_1$, we have (a)
$x$ and $y$ are at distance one in $I_1$ if and only if $f(x)$ and
$f(y)$ are at distance two in $I_2$, and (b) $x$ and $y$ are at
distance two in $I_1$ if and only if $f(x)$ and $f(y)$ are at
distance one in $I_2$. (It follows that $x$ and $y$ are at
distance 3 (antipodal) in $I_1$ if and only if $f(x)$ and $f(y)$
are at distance 3 (antipodal) in $I_2$.) Here distance refers to
the usual graphical distance on the respective edge graph. In case
$V(I_1) = V(I_2)$ and the identity map is an antimorphism between
$I_1$ and $I_2$, then we say that $I_1$ and $I_2$ are {\em
antimorphic}. Thus, the two icosahedra in Figure 2 below are
antimorphic (the map, taking each vertex of the left icosahedron
in Figure 2 to the vertex of the same name in the right
icosahedron, is an antimorphism).}
\end{defn}

\setlength{\unitlength}{5mm}

\begin{picture}(29,15.5)(1,-0.5)

%%%%%%%%%%%%%%%%%%%%%%% I_{1} %%%%%%%%%%%%%%%%

\thicklines

\put(7,2){\line(-3,2){6}} \put(7,2){\line(3,2){6}}
\put(7,2){\line(-2,3){2}} \put(7,2){\line(2,3){2}}

\put(1,6){\line(4,-1){4}} \put(1,6){\line(0,1){4}}
\put(5,5){\line(1,0){4}} \put(5,5){\line(-4,5){4}}
\put(5,5){\line(1,2){2}}

\put(9,5){\line(4,1){4}} \put(9,5){\line(-1,2){2}}
\put(9,5){\line(4,5){4}} \put(13,6){\line(0,1){4}}

\put(7,9){\line(-6,1){6}} \put(7,9){\line(6,1){6}}
\put(7,9){\line(0,1){5}}

\put(7,14){\line(-3,-2){6}} \put(7,14){\line(3,-2){6}}

\thinlines

\put(0,-0.5){\line(1,0){30}} \put(0,15){\line(1,0){30}}
%%These two are boundary lines

\put(5,11){\line(-4,-1){4}} \put(5,11){\line(2,3){2}}
\put(5,11){\line(1,0){1.7}} \put(9,11){\line(-1,0){1.7}}

\put(9,11){\line(4,-1){4}} \put(9,11){\line(-2,3){2}}

\put(1,6){\line(4,5){1.4}} \put(2.8,8.25){\line(4,5){0.9}}
\put(5,11){\line(-4,-5){0.938}}

\put(13,6){\line(-4,5){1.4}} \put(11.2,8.25){\line(-4,5){0.9}}
\put(9,11){\line(4,-5){0.938}}

\put(7,2){\line(0,1){2.7}} \put(7,7){\line(0,-1){1.7}}

\put(1,6){\line(6,1){2.5}} \put(4.3,6.55){\line(6,1){1.18}}

\put(13,6){\line(-6,1){2.5}} \put(9.7,6.55){\line(-6,1){1.18}}

\put(5,11){\line(1,-2){0.75}} \put(9,11){\line(-1,-2){0.75}}

\put(7,7){\line(-6,-1){0.8}} \put(7,7){\line(6,-1){0.8}}

\put(7,7){\line(-1,2){0.4}} \put(6.4,8.2){\line(-1,2){0.36}}

\put(7,7){\line(1,2){0.4}} \put(7.6,8.2){\line(1,2){0.36}}

\put(7.6,14){$x_{23}$} \put(7.6,1.5){$x_{41}$}

\put(0,10.5){$x_{12}$} \put(7.2,9.5){$x_{14}$}
\put(13,10.5){$x_{21}$} \put(9.1,11.3){$x_{24}$}
\put(4,11.3){$x_{31}$}

\put(0,5.5){$x_{34}$} \put(4,4.4){$x_{13}$}
\put(9.1,4.4){$x_{42}$} \put(13.1,5.5){$x_{43}$}
\put(5.9,6.2){$x_{32}$}

\put(3.5,1.5){\mbox{$I_{1}$}}

%%%%%%%%%%%%%%%% I_{2} %%%%%%%%%%%%%%%%%%%%%

\thicklines

\put(22,2){\line(-3,2){6}} \put(22,2){\line(3,2){6}}
\put(22,2){\line(-2,3){2}} \put(22,2){\line(2,3){2}}

\put(16,6){\line(4,-1){4}} \put(16,6){\line(0,1){4}}
\put(20,5){\line(1,0){4}} \put(20,5){\line(-4,5){4}}
\put(20,5){\line(1,2){2}}

\put(24,5){\line(4,1){4}} \put(24,5){\line(-1,2){2}}
\put(24,5){\line(4,5){4}} \put(28,6){\line(0,1){4}}

\put(22,9){\line(-6,1){6}} \put(22,9){\line(6,1){6}}
\put(22,9){\line(0,1){5}}

\put(22,14){\line(-3,-2){6}} \put(22,14){\line(3,-2){6}}

\thinlines

\put(20,11){\line(-4,-1){4}} \put(20,11){\line(2,3){2}}
\put(20,11){\line(1,0){1.7}} \put(24,11){\line(-1,0){1.7}}

\put(24,11){\line(4,-1){4}} \put(24,11){\line(-2,3){2}}

\put(16,6){\line(4,5){1.4}} \put(17.8,8.25){\line(4,5){0.9}}
\put(20,11){\line(-4,-5){0.938}}

\put(28,6){\line(-4,5){1.4}} \put(26.2,8.25){\line(-4,5){0.9}}
\put(24,11){\line(4,-5){0.938}}

\put(22,2){\line(0,1){2.7}} \put(22,7){\line(0,-1){1.7}}

\put(16,6){\line(6,1){2.5}} \put(19.3,6.55){\line(6,1){1.18}}

\put(28,6){\line(-6,1){2.5}} \put(24.7,6.55){\line(-6,1){1.18}}

\put(20,11){\line(1,-2){0.75}} \put(24,11){\line(-1,-2){0.75}}

\put(22,7){\line(-6,-1){0.8}} \put(22,7){\line(6,-1){0.8}}

\put(22,7){\line(-1,2){0.4}} \put(21.4,8.2){\line(-1,2){0.36}}

\put(22,7){\line(1,2){0.4}} \put(22.6,8.2){\line(1,2){0.36}}

\put(22.6,14){$x_{12}$} \put(22.6,1.5){$x_{43}$}

\put(15,10.5){$x_{21}$} \put(22.2,9.5){$x_{41}$}
\put(28,10.5){$x_{24}$} \put(24.1,11.3){$x_{42}$}
\put(19,11.3){$x_{32}$}

\put(15,5.5){$x_{13}$} \put(19,4.4){$x_{31}$}
\put(24.1,4.4){$x_{14}$} \put(28.1,5.5){$x_{34}$}
\put(20.9,6.2){$x_{23}$}

\put(18.5,1.5){\mbox{$I_{2}$}}

\put(5.7,0){\mbox{\bf Figure 2\,: An antimorphic pair of
icosahedra}}

\end{picture}

\medskip

\noindent {\bf Another description of {\boldmath $(S^{\,2} \times
S^{\,2})_{12}$}\,:} Take an antimorphic pair of icosahedra, say
$I_1$ and $I_2$ (with common vertex set $V$). It turns out that
$I_1$ and $I_2$ have the identical automorphism group $A_5 \times
\ZZ_2$ (not merely isomorphic, cf. Lemma \ref{L3} below). Also,
there is a bijection $\varphi$ from the triangles of $I_1$ to the
triangles of $I_2$ such that for each triangle $\Delta = abc$ of
$I_1$, $\varphi(\Delta)= ijk$ is the only triangle of $I_2$ for
which $aij$, $bjk$ and $cik$ are triangles of $I_2$ (cf. Lemma
\ref{L3}). Now, the vertex-set of $(S^{\,2} \times S^{\,2})_{12}$
is $V$ ($= V(I_1) = V(I_2)$) and it has two types of facets. (i)
For each vertex $x$, the neighbors of $x$ in $I_1$ form facets.
(ii) For each triangle $\Delta$ of $I_1$ and each vertex $y$ in
$\Delta^{\prime} = \varphi(\Delta)$, $(\Delta \cup
\Delta^{\prime}) \setminus \{y\}$ is a facet. Thus $(S^{\,2}
\times S^{\,2})_{12}$ has 12 facets of the first type and $20
\times 3 = 60$ facets of the second type. From the description,
it is clear that the
common automorphism group $A_5 \times \ZZ_2$ of $I_1$ and $I_2$ is
an automorphism group of $(S^{\,2} \times S^{\,2})_{12}$. It turns
out that its full automorphism group is $2S_5$
generated by the two automorphisms $g = (x_{12} x_{21} x_{24}
x_{42} x_{14} x_{41} x_{43} x_{34} x_{13} x_{31} x_{32} x_{23})$
and $h = (x_{12} x_{14} x_{21} x_{24} x_{31}) (x_{13} x_{42}
x_{43} x_{32} x_{34})$. The automorphism $g$ interchanges $I_1$
and $I_2$.

\begin{remark}$\!\!\!${\bf .} \label{R4}
{\rm It should be emphasized that the existence of an antimorphic
pair of icosahedra (exploited in the above construction of $(S^{\,2}
\times S^{\,2})_{12}$) is a minor miracle, and only an empirically
verified fact. Its deeper geometric significance, if any, remains
to be understood. }
\end{remark}

\section{A self-dual CW decomposition of {\boldmath $\CCC
P^{\hspace{.1mm}2}$}}

Here we have taken the cell complex
$\partial T_0 \times \partial T_0$, and triangulated it to obtain
the simplicial complex $(S^{\,2}\times S^{\,2})_{16}$ and finally
quotiented this simplicial complex by $\ZZ_2$ to obtain $\CC
P^{\,2}_{10}$. This procedure reflects our obsession with
simplicial complexes. However, one may straightaway quotient the
cell complex by $\ZZ_2$ to obtain a (non-regular) CW decomposition
of $\CC P^{\,2}$. This CW complex is self-dual in the sense that
its face-vector $(10, 24, 31, 24, 10)$ exhibits a curious
palindromic symmetry. We proceed to describe it in some details.
Consider the $\ZZ_2$ action on $\RR^{\,6} \equiv \RR^{\,3} \times
\RR^{\,3}$ given by $(x, y) \leftrightarrow (y, x)$. Let $\eta
\colon \RR^{\,6} \to \RR^{\,6}/\ZZ_2$ be the quotient map. We know
that $\eta(S^{\,2}\times S^{\,2}) = \CC P^{\,2}$. We give a CW
decomposition $W$ of the space $\eta(\partial T_0 \times \partial
T_0)$.

For $0\leq i\leq 4$, let $W^i$ denote the set of $i$-cells in $W$.
For $i\neq 2$ the $i$-cells in $W$ are the images (under the map
$\eta$) of $i$-cells in $\partial T_0 \times \partial T_0$. A
2-cell in $W$ is the image of a 2-cell $F$ in $\partial T_0 \times
\partial T_0$ which is not of the form $E \times E$ for some edge
$E$ in $\partial T_0$. More explicitly
\begin{eqnarray*}
W^0 & = & V(\CC P^{\,2}_{10}), \\
W^1 & = & \{\eta(E) \, : \, E \mbox{ is an edge of } \,
\partial T_0 \times \partial T_0\}, \\
W^2 & = & \{\eta(|x_{ij}x_{ik}x_{il}|) \, : \, 1\leq j < k < l
\leq 4, \, 1\leq i \leq 4\}   \\
&& ~~ \cup \{\eta(|x_ix_j| \times |x_kx_l|) \, : \, i<j, \, k<l \,
\mbox{ and  either } i< k \mbox{ or } i=k \mbox{ and } j < l\}, \\
W^3 & = &  \{\eta(A) \, : \, A \mbox{ is a 3-cell of } \,
\partial T_0 \times \partial T_0\} \mbox{ and} \\
 W^4 & = & \{\eta(B) \, : \, B \mbox{ is a 4-cell of } \,
\partial T_0 \times \partial T_0\}.
\end{eqnarray*}

Then, $W^1$ contains 24 cells, $W^2$ contains $16+15 = 31$ cells,
$W^3$ contains $4\times 6 =24$ cells and $W^4$ contains 10 cells.
Clearly, each 1-cell in $W$ is regular (i.e., homeomorphic to a
closed interval). Since all the 2-cells are homeomorphic images of
the corresponding 2-cells in $\partial T_0 \times \partial T_0$,
it follows that all the 2-cells in $W$ are regular.

For $0\leq i\leq 4$, let $X_i = \bigcup_{\beta \in W^0 \cup \cdots
\cup W^i} \beta$. Then $\partial \alpha \subseteq X_{i-1}$ if
$\alpha \in W^i$ for $i \neq 3$. Let $\gamma$ be a 3-cell in $W$.
If $\gamma = \eta(|x_ix_jx_k| \times |x_ix_j|)$, $i<j<k$, then
$\gamma$ is obtained from $|x_i x_jx_k| \times |x_ix_j|$ by
identifying $|x_{ii}x_{jj}x_{ij}|$ with $|x_{ii}x_{jj} x_{ji}|$
(by the identification given by $x_{ij} \leftrightarrow x_{ji}$).
Thus, $\gamma$ is a regular 3-cell and $\partial \gamma =
\eta(|x_i x_k| \times |x_ix_j|) \cup \eta(|x_jx_k| \times
|x_ix_j|) \cup \eta(|x_{ii}x_{ji}x_{ki}|) \cup \eta(|x_{ij}x_{jj}
x_{kj}|)$. (Now, it is clear why we do not have to take
$\eta(|x_ix_j| \times |x_ix_j|)$ in $W^2$. In fact, $\eta(|x_ix_j|
\times |x_ix_j|)$ is inside of $\gamma$.) Therefore, $\partial
\gamma \subseteq X_2$. Same things are true if $\gamma =
\eta(|x_ix_jx_k| \times |x_ix_k|)$ or $\eta(|x_ix_jx_k| \times
|x_jx_k|)$. On the other hand, if $\gamma = \eta(F \times E)$,
where $E$ is an edge and $F$ is a 2-simplex and $E \not\subseteq
F$, then $\gamma$ is homeomorphic to $F \times E$ and hence is a
regular 3-cell. In this case, it follows from the definition of
$W^2$ that $\partial \gamma \subseteq X_2$. Thus $W$ is a CW
complex.

If $\sigma$ is a 4-cell in $W$ then, either $\sigma = \eta(|x_ix_j
x_k| \times |x_ix_jx_k|)$, for some $i<j<k$ or $\sigma = \eta(|x_i
x_jx_k| \times |x_ix_j x_l|)$, where $\{i, j, k, l\}$ is an even
permutation of $\{1, 2, 3, 4\}$. In the first case, $\sigma$ is
homeomorphic to $|x_{ii}x_{jj}x_{kk} x_{ij}x_{ik}| \cup |x_{ii}
x_{jj}x_{kk}x_{ij}x_{jk}| \cup |x_{ii}x_{jj}x_{kk} x_{ik}x_{jk}|$
and hence $\sigma$ is a regular 4-cell. In the second case,
$\sigma$ is obtained from $|x_ix_jx_k| \times |x_ix_j x_l|$ by
identifying $|x_{ii}x_{jj}x_{ij}|$ with $|x_{ii}x_{jj} x_{ji}|$
(by the identification given by $x_{ij} \leftrightarrow x_{ji}$).
So, $\sigma$ is not a regular cell. Thus $W^4$ contains four
regular 4-cells and six singular 4-cells.

Since each cell in $W$ is the quotient of a cell in $S^{\,2}_4
\times S^{\,2}_4$, $(S^{\,2} \times S^{\,2})_{16}$ is a simplicial
subdivision of $S^{\,2}_4 \times S^{\,2}_4$ and $\CC P^{\,2}_{10}$
is the quotient of $(S^{\,2} \times S^{\,2})_{16}$, it follows that
$\CC P^{\,2}_{10}$ is a simplicial subdivision of $W$.

\section{Proofs}

\begin{defn}$\!\!\!${\bf .} \label{pure}
{\rm Let $G$ be a group of simplicial automorphisms of a
simplicial complex $X$ with vertex set $V(X)$. We shall say that the
action of $G$ on $X$ is {\em pure} if it satisfies\,: $(a)$
whenever $u$, $v$ are distinct vertices from the same $G$-orbit,
$uv$ is a non-edge of $X$, and $(b)$ for each $G$-orbit $\theta
\subseteq V(X)$ and each $\alpha \in X$, the stabiliser
$G_{\alpha}$ of $\alpha$ in $G$ acts transitively on $\theta \cap
V({\rm lk}_X(\alpha))$.}
\end{defn}

\begin{lemma}$\!\!\!${\bf .} \label{L1}
Let $G$ be a group of simplicial automorphisms of a simplicial
complex $X$. Let $q \colon V(X) \to V(X)/G$ denote the quotient
map, and $X/G := \{q(\alpha) \, : \, \alpha\in X\}$. If the action
of $G$ on $X$ is pure then $X/G$ is a simplicial complex which
triangulates $|X|/G$ $($where the action of $\,G$ on $V(X)$ is
extended to an action of $\,G$ on $|X|$ piecewise linearly, i.e.,
affinely on the geometric carrier of each simplex$)$. That is,
we have $|X/G| = |X|/G$.
\end{lemma}

\noindent {\bf Proof.} The condition $(a)$ ensures that the
quotient map $q$ is one-one on each simplex of $X$. The simplicial
map $q \colon X \to X/G$ induces a piecewise linear continuous map
$|q|$ from $|X|$ onto $|X/G|$.

\smallskip

\noindent {\sf Claim.} {\em The fibres of $\,q \colon X \to X/G$
are precisely the $G$-orbits on simplices of $X$ $($that is, if
$\alpha$, $\alpha^{\,\prime} \in X$ are such that $q(\alpha) =
q(\alpha^{\,\prime})$ then there exists $g \in G$ such that
$g(\alpha) = \alpha^{\,\prime})$.}

\smallskip

We prove the claim by induction on $k = \dim(\alpha) =
\dim(\alpha^{\,\prime})$. The claim is trivial for $k = -1$. So,
assume $k \geq 0$, and the claim is true for all smaller
dimensions. Choose a simplex $\beta \subseteq \alpha$ of dimension
$k-1$, and let $\beta^{\,\prime}\subseteq \alpha^{\,\prime}$ be
such that $q(\beta^{\,\prime}) = q(\beta)$. By induction
hypothesis, $\beta^{\,\prime}$ and $\beta$ are in the same
$G$-orbit. Therefore, applying a suitable element of $G$, we may
assume, without loss of generality, that $\beta^{\,\prime} =
\beta$. Let $\alpha = \beta \cup \{x\}$, $\alpha^{\,\prime} =
\beta \cup \{x^{\,\prime}\}$. Then $q(x) = q(x^{\,\prime})$, i.e.,
$x$ and $x^{\,\prime}$ are in the same $G$-orbit. Now, by
assumption $(b)$, there is a $g\in G_{\beta}$ such that $g(x)=
x^{\,\prime}$. Then $g(\alpha) =\alpha^{\,\prime}$. This proves
the claim.

\smallskip

In the presence of condition $(a)$, the claim ensures that the
fibres of $|q|$ are precisely the $G$-orbits on points of $|X|$.
Hence $|q|$ induces the required homeomorphism between $|X|/G$ and
$|X/G|$. \hfill $\Box$

\medskip

Up to isomorphism, there are exactly two 6-vertex 2-spheres,
namely, $S_1$ and $S_2$ given in Figure 3. We need the following
lemma to prove Theorem \ref{T1}.

\bigskip

\setlength{\unitlength}{4mm}

\begin{picture}(36,13)(0,1)

%%%%%%%%%%%%%%%%%%%%%%%%%%%%%%%%%%%%%%%%%%%%%%%%%%%%
\thicklines

\put(2,5){\line(1,0){12}} \put(2,5){\line(3,4){6}}
\put(2,5){\line(2,1){4}} \put(2,5){\line(6,5){6}}
\put(6,7){\line(1,0){4}} \put(6,7){\line(2,3){2}}
\put(10,7){\line(-1,3){2}} \put(10,7){\line(-2,3){2}}
\put(8,10){\line(0,1){3}} \put(14,5){\line(-3,4){6}}
\put(14,5){\line(-4,1){8}} \put(14,5){\line(-2,1){4}}

\put(5.5,6){\mbox{$a_{1}$}} \put(10.3,7.2){\mbox{$a_{2}$}}
\put(6.8,10){\mbox{$a_{3}$}} \put(1.2,5.5){\mbox{$b_{1}$}}
\put(14,5.5){\mbox{$b_{2}$}} \put(6.5,12.5){\mbox{$b_{3}$}}

\put(4,3.5){\mbox{$S_1= S^0_2\ast S^0_2\ast S^0_2$}}

%%%%%%%%%%%%%%%%%%%%%%%%%%%%%%%%%%%%%%%%%%%%%%%%%%%%
\thicklines

\put(18,5){\line(-3,4){6}} \put(18,5){\line(3,4){6}}
\put(18,5){\line(-1,3){2}} \put(18,5){\line(1,3){2}}
\put(18,5){\line(0,1){3}}

\put(18,8){\line(-2,3){2}} \put(18,8){\line(2,3){2}}
\put(16,11){\line(-2,1){4}} \put(16,11){\line(4,1){8}}
\put(16,11){\line(1,0){4}} \put(12,13){\line(1,0){12}}
\put(20,11){\line(2,1){4}}

\put(11.2,12){\mbox{$b_{1}$}} \put(24.2,12){\mbox{$b_{2}$}}
\put(16.5,5){\mbox{$b_{3}$}} \put(15,10.3){\mbox{$a_{1}$}}
\put(20.3,10.3){\mbox{$a_{2}$}} \put(17.6,9.1){\mbox{$a_{3}$}}

\put(19.5,4.5){\mbox{$S_2$}}

%%%%%
%%%%%%%%%%%%%%%%%%%%%%%%%%%%%%%%%%%%%%%%%%%%%%%%%%%%
\thicklines

\put(32,5){\line(-3,1){3}} \put(32,5){\line(3,1){3}}
\put(32,5){\line(0,1){5}} \put(29,6){\line(1,0){1}}
\put(30.5,6){\line(1,0){1}} \put(34,6){\line(1,0){1}}
\put(32.2,6){\line(1,0){1}} \put(29,6){\line(0,1){5}}
\put(35,6){\line(0,1){5}}

\put(32,10){\line(-3,1){3}} \put(32,10){\line(3,1){3}}
\put(29,11){\line(1,0){6}}

\put(32,4){\mbox{$b_{3}$}} \put(32.3,9.4){\mbox{$a_{3}$}}
\put(28,6){\mbox{$b_{1}$}} \put(27.9,10.5){\mbox{$a_{1}$}}
\put(35.4,6){\mbox{$b_{2}$}} \put(35.4,10.5){\mbox{$a_{2}$}}

\put(30,3.5){\mbox{$C$}}

\put(1,1.7){\mbox{\bf Figure 3\,(a): 6-vertex 2-spheres}}

\put(21,1.7){\mbox{\bf Figure 3\,(b): Triangular prism}}

\thinlines

\put(-1,1){\line(1,0){37}} \put(-1,14){\line(1,0){37}}
%%These two are boundary lines

\end{picture}

\begin{lemma}$\!\!\!${\bf .} \label{L2}
Let $C$ be the triangular prism given in Figure $3$\,{\rm (b)} $($i.e., $C$
is the product of a $2$-simplex and an edge$)$. Up to isomorphism,
there exists a unique $6$-vertex simplicial subdivision
$\widehat{C}$ of $C$. The facets $($tetrahedra$)$ in
$\widehat{C}$ are $a_1b_1b_2 b_3, a_1a_2b_2b_3, a_1a_2a_3b_3$.
Moreover, $\partial \widehat{C}$ is isomorphic to $S_2$ of
Figure $3$\,{\rm (a)} and determines $\widehat{C}$ uniquely.
\end{lemma}

\noindent {\bf Proof.} Let $\widehat{C}$ be a 6-vertex
subdivision of $C$. Then there exists a 3-simplex $\sigma$ in
$\widehat{C}$ which contains the 2-simplex $b_1b_2b_3$. Without
loss of generality, we may assume that $\sigma = a_1b_1b_2b_3$.
Then $C$ is the union of $\sigma$ and the pyramid $P$ given in
Figure 4. Since we are not allowed to introduce new vertices,
clearly the rectangular base of $P$ must be triangulated using two
triangles, in one of two isomorphic ways, and the remaining
tetrahedra in $\widehat{C}$ must have the apex of $P$ as a
vertex and one of these two triangles as base. Thus, without loss
of generality, $P = a_1a_2b_2b_3 \cup a_1a_2a_3b_3$. This proves
the first part.

\bigskip

\setlength{\unitlength}{4mm}

\begin{picture}(36.5,12)(3,0.5)

%%%%%%%%%%%%%%%%%%%%%%%%%%%%%%%%%%%%%%%%%%%%%%%%%%%%
\thicklines

\put(6,5){\line(-3,1){3}} \put(6,5){\line(3,1){3}}
\put(3,6){\line(1,0){0.75}} \put(4.25,6){\line(1,0){0.75}}
\put(7.75,6){\line(1,0){1.25}} \put(6,6){\line(1,0){1.25}}
\put(3,6){\line(0,1){5}}

\put(11,5){\line(3,1){3}} \put(11,5){\line(0,1){5}}
\put(14,6){\line(0,1){5}} \put(11,10){\line(-3,1){3}}
\put(11,10){\line(3,1){3}} \put(8,11){\line(1,0){6}}

\thinlines

\put(2,0.5){\line(1,0){37}} \put(2,12.5){\line(1,0){37}}
%%These two are boundary lines

\put(3,11){\line(6,-5){6}} \put(6,5){\line(-1,2){3}}
\put(11,5){\line(-1,2){3}}

\put(8,11){\line(6,-5){1.1}} \put(9.5,9.75){\line(6,-5){1.1}}
\put(11.3,8.25){\line(6,-5){0.88}} \put(14,6){\line(-6,5){1.5}}

\put(5.5,4.2){\mbox{$b_{3}$}} \put(3.4,6.5){\mbox{$b_{1}$}}
\put(3.5,11){\mbox{$a_{1}$}} \put(8.7,5){\mbox{$b_{2}$}}

\put(10,4.5){\mbox{$b_{3}$}} \put(11.3,9.4){\mbox{$a_{3}$}}
\put(8,11.5){\mbox{$a_{1}$}} \put(13.5,5){\mbox{$b_{2}$}}
\put(13.5,11.5){\mbox{$a_{2}$}} \put(7.9,8){\mbox{$\cup$}}

\put(5.5,2.5){\mbox{$C = a_1b_1b_2b_3 \cup P$}}

%%%%%%%%%%%%%%%%%%%%%%%%%%%%%%%%%%%%%%%%%%%%%%

\thicklines

\put(19,5){\line(0,1){5}} \put(19,10){\line(-3,1){3}}
\put(19,10){\line(3,1){3}} \put(16,11){\line(1,0){6}}

\put(26,5){\line(3,1){3}} \put(29,6){\line(0,1){5}}
 \put(23,11){\line(1,0){6}}

\thinlines

\put(19,5){\line(-1,2){3}} \put(19,5){\line(1,2){3}}
\put(26,5){\line(-1,2){3}} \put(26,5){\line(1,2){3}}

\put(23,11){\line(6,-5){1.1}} \put(24.5,9.75){\line(6,-5){1.1}}
\put(26,8.5){\line(6,-5){0.88}} \put(29,6){\line(-6,5){1.5}}

\put(16,11.5){\mbox{$a_{1}$}} \put(23,11.5){\mbox{$a_{1}$}}
\put(21.5,11.5){\mbox{$a_{2}$}} \put(28.5,11.5){\mbox{$a_{2}$}}
\put(19.3,9.4){\mbox{$a_{3}$}} \put(18,4.5){\mbox{$b_{3}$}}
\put(25,4.5){\mbox{$b_{3}$}} \put(28.5,5){\mbox{$b_{2}$}}
\put(22,7){\mbox{$\cup$}}

\put(18,2.8){\mbox{$P = a_1a_2a_3b_3 \cup a_1a_2b_2b_3$}}

%%%%%%%%%%%%%%%%%%%%%%%%%%%%%%%%%%%%%%%%%%%%%%%%%%%%
\thicklines

\put(34,5){\line(-3,1){3}} \put(34,5){\line(3,1){3}}
\put(34,5){\line(0,1){5}} \put(31,6){\line(1,0){0.75}}
\put(32.25,6){\line(1,0){0.75}} \put(36.25,6){\line(1,0){0.75}}
\put(35,6){\line(1,0){0.75}} \put(31,6){\line(0,1){5}}
\put(37,6){\line(0,1){5}}

\put(34,10){\line(-3,1){3}} \put(34,10){\line(3,1){3}}
\put(31,11){\line(1,0){6}}

\thinlines

\put(34,5){\line(1,2){3}} \put(34,5){\line(-1,2){3}}

\put(37,6){\line(-6,5){1.5}} \put(31,11){\line(6,-5){1.1}}
\put(32.5,9.75){\line(6,-5){1.1}}
\put(34.15,8.375){\line(6,-5){0.88}} %(32.3,8.25)

\put(33,4.2){\mbox{$b_{3}$}} \put(34.3,9.4){\mbox{$a_{3}$}}
\put(31.4,6.5){\mbox{$b_{1}$}} \put(31,11.5){\mbox{$a_{1}$}}
\put(36.5,5){\mbox{$b_{2}$}} \put(36.5,11.5){\mbox{$a_{2}$}}

\put(35,3){\mbox{$\widehat{C}$}}

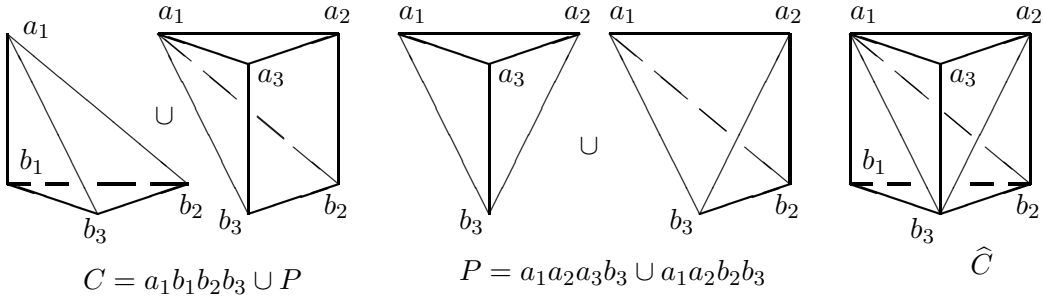
\put(7.5,1){\mbox{\bf Figure 4\,: Simplicial subdivision of the
triangular prism}}
\end{picture}

\medskip

The last part follows from the fact that the facets of
$\widehat{C}$ are the maximal cliques in the 1-skeleton of
$\partial \widehat{C}$. \hfill $\Box$

\bigskip

\noindent {\bf Proof of Theorem \ref{T1}.} Let $X$ be a
$16$-vertex simplicial subdivision of $S^{\,2}_4 \times S^{\,2}_4$
satisfying $(i)$, $(ii)$ and $(iii)$.

For $i\neq j$, consider the 2-cell $x_ix_j \times x_ix_j$. By
$(iii)$, $x_{ij} x_{ji}$ can not be an edge in
$X$. This implies that $x_{ii}x_{jj}$, $x_{ii}x_{ji}x_{jj}$,
$x_{ii} x_{ij} x_{jj} \in X$ and $x_ix_j \times x_ix_j =
x_{ii}x_{ji} x_{jj} \cup x_{ii}x_{ij}x_{jj}$ (cf. Figure 5\,$(a)$).

For $i, j, k$ distinct, consider the 2-cell $x_ix_j\times x_ix_k$.
Since $X$ satisfies $(iii)$, both $x_{ij}$ and
$x_{ji}$ can't be in ${\rm lk}_X(x_{ik})$. Now, $x_{ik}x_{ij}$ is
an edge in the cell complex $S^{\,2}_4 \times S^{\,2}_4$ and hence
is an edge in $X$. Thus, $x_{ik} x_{ji}$ can not be an edge in
$X$. This implies that $x_{ii}x_{jk}$, $x_{ii}x_{ji} x_{jk}$,
$x_{ii}x_{ik}x_{jk} \in X$ and $x_ix_j \times x_ix_k =
x_{ii}x_{ji}x_{jk} \cup x_{ii} x_{ik}x_{jk}$ (cf. Figure 5\,$(b)$).

\bigskip

\setlength{\unitlength}{3mm}

\begin{picture}(48.5,14)(2,-3.5)

%%%%%%%%%%%%%%%%%%%%%%%%%%%%%%%%%%%%%%%%%%%%%%%%%%%%
\thicklines

\put(3,4){\line(1,0){5}} \put(3,4){\line(0,1){5}}
\put(3,9){\line(1,0){5}} \put(8,4){\line(0,1){5}}

\thinlines

\put(1,-3.5){\line(1,0){49}} \put(1,11){\line(1,0){49}}
%%These two are boundary lines

\put(3,4){\line(1,1){5}}

\put(2,3){\mbox{$x_{ii}$}} \put(7.4,3){\mbox{$x_{ji}$}}
\put(2,9.7){\mbox{$x_{ij}$}} \put(7.4,9.7){\mbox{$x_{jj}$}}
\put(2.2,1.3){\mbox{$x_{i}x_j\times x_ix_j$}}

%%%%%

\thicklines

\put(13,4){\line(1,0){5}} \put(13,4){\line(0,1){5}}
\put(13,9){\line(1,0){5}} \put(18,4){\line(0,1){5}}

\thinlines

\put(13,4){\line(1,1){5}}

\put(12,3){\mbox{$x_{ii}$}} \put(17.4,3){\mbox{$x_{ji}$}}
\put(12,9.7){\mbox{$x_{ik}$}} \put(17.4,9.7){\mbox{$x_{jk}$}}
\put(12.2,1.3){\mbox{$x_ix_j\times x_ix_k$}}

%%%%%%

\thicklines

\put(23,4){\line(1,0){5}} \put(23,4){\line(0,1){5}}
\put(23,9){\line(1,0){5}} \put(28,4){\line(0,1){5}}

\thinlines

\put(23,4){\line(1,1){5}}

\put(22,3){\mbox{$x_{12}$}} \put(27.4,3){\mbox{$x_{32}$}}
\put(22,9.7){\mbox{$x_{14}$}} \put(27.4,9.7){\mbox{$x_{34}$}}
\put(22.2,1.3){\mbox{$x_{1}x_3 \times x_2x_4$}}

%%%%%

\thicklines

\put(33,4){\line(1,0){5}} \put(33,4){\line(0,1){5}}
\put(33,9){\line(1,0){5}} \put(38,4){\line(0,1){5}}

\thinlines

\put(38,4){\line(-1,1){5}}

\put(32,3){\mbox{$x_{21}$}} \put(37.4,3){\mbox{$x_{31}$}}
\put(32,9.7){\mbox{$x_{24}$}} \put(37.4,9.7){\mbox{$x_{34}$}}
\put(32.2,1.3){\mbox{$x_{2}x_3\times x_1x_4$}}

%%%%%

\thicklines

\put(43,4){\line(1,0){5}} \put(43,4){\line(0,1){5}}
\put(43,9){\line(1,0){5}} \put(48,4){\line(0,1){5}}

\thinlines

\put(48,4){\line(-1,1){5}}

\put(42,3){\mbox{$x_{13}$}} \put(47.4,3){\mbox{$x_{23}$}}
\put(42,9.7){\mbox{$x_{14}$}} \put(47.4,9.7){\mbox{$x_{24}$}}
\put(42.2,1.3){\mbox{$x_{1}x_2\times x_3x_4$}}

%%%%%

\put(4,-0.6){\boldmath $(a)$} \put(14,-0.6){\boldmath $(b)$}
\put(24,-0.6){\boldmath $(c)$} \put(34,-0.6){\boldmath $(d)$}
\put(44,-0.6){\boldmath $(e)$}

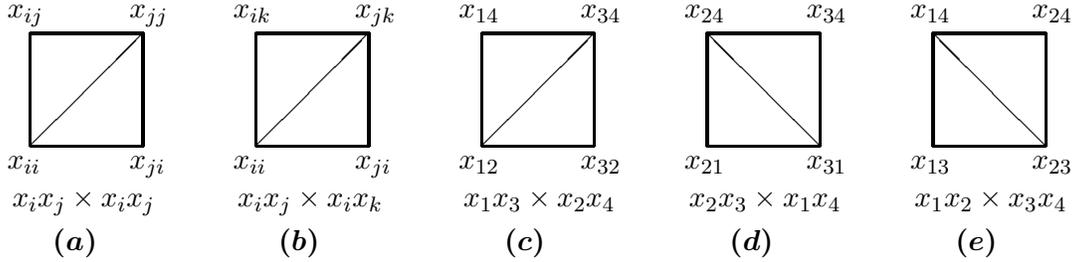
\put(2,-2.5){\mbox{\bf Figure 5\,: Simplicial subdivisions of
rectangular 2-cells of {\boldmath $S^{\,2}_4 \times S^{\,2}_4$}}}

\end{picture}

\medskip

Consider the 2-cell $x_1x_3\times x_2x_4$. Clearly, $x_1x_3 \times
x_2x_4 = x_{12}x_{32}x_{34} \cup x_{12}x_{14}x_{34}$ or $=
x_{12}x_{32}x_{14} \cup x_{32}x_{14}x_{34}$.

\smallskip

\noindent {\bf Case 1.} $x_1x_3\times x_2x_4 = x_{12}x_{32}x_{34}
\cup x_{12}x_{14}x_{34}$ (cf. Figure 5\,$(c)$). So, $x_{12}x_{34}\in X$. Then, by (ii),
$x_{21}x_{43}\in X$ and, by (iii), $x_{12}x_{43}$, $x_{21}x_{34}
\not\in X$. This implies that $x_2x_3\times x_1x_4=
x_{21}x_{31}x_{24} \cup x_{31}x_{24}x_{34}$ (cf. Figure 5\,$(d)$). So, $x_{31}x_{24}\in
X$. Then, by (ii), $x_{13}x_{42}\in X$ and, by (iii),
$x_{13}x_{24}$, $x_{31}x_{42} \not\in X$. This implies that
$x_1x_2\times x_3x_4= x_{13}x_{23}x_{14} \cup x_{23}x_{14}x_{24}$ (cf. Figure 5\,$(e)$).
So, $x_{14}x_{23} \in X$. Then, by (ii), $x_{41}x_{32}\in X$ and,
by (iii), $x_{14}x_{32}$, $x_{41}x_{23} \not\in X$. These give the
2-skeleton of $X$. Observe that we have already 84 edges as
mentioned in the construction of $(S^{\,2} \times S^{\,2})_{16}$
and, since $X$ satisfies $(iii)$, all the 36
remaining 2-sets are non-edges in $X$.

Observe that any 3-cell in $S^{\,2}_4 \times S^{\,2}_4$ is the
product of a 2-simplex and an edge. For $i, j, k$ distinct,
consider the 3-cell $x_ix_jx_k \times x_ix_j$. Since
$x_{ii}x_{jj}$, $x_{ii}x_{ki}$ and $x_{jj}x_{ki}$ are edges, by
Lemma \ref{L2}, $x_ix_jx_k \times x_ix_j =
x_{ii}x_{ij}x_{kj}x_{jj} \cup x_{ii}x_{ki}x_{kj}x_{jj} \cup
x_{ii}x_{ki}x_{ji}x_{jj}$ is the unique subdivision of $x_ix_jx_k
\times x_ix_j$ (cf. Figure 6\,$(a)$). Similarly, $x_ix_j \times x_ix_jx_k =
x_{ii}x_{ji}x_{jk}x_{jj}  \cup x_{ii}x_{ik}x_{jk}x_{jj} \cup
x_{ii}x_{ik}x_{ij}x_{jj}$ is the unique subdivision of $x_ix_j
\times x_ix_jx_k$ (cf. Figure 6\,$(b)$).

\medskip

\setlength{\unitlength}{4mm}

\begin{picture}(36.5,14)(3,-1)

%%%%%%%%%%%%%%%%%%%%%%%%%%%%%%%%%%%%%%%%%%%%%%%%%%%%

%%%%%%%%%%%%%%%%%%%%%%%%%%%%%%%%%%%%%%%%%%%%%%%%%%%%
\thicklines

\put(6,5){\line(-3,1){3}} \put(6,5){\line(3,1){3}}
\put(6,5){\line(0,1){5}} \put(3,6){\line(1,0){0.75}}
\put(4.25,6){\line(1,0){0.75}} \put(8.25,6){\line(1,0){0.75}}
\put(7,6){\line(1,0){0.75}} \put(3,6){\line(0,1){5}}
\put(9,6){\line(0,1){5}}

\put(6,10){\line(-3,1){3}} \put(6,10){\line(3,1){3}}
\put(3,11){\line(1,0){6}}

\thinlines

\put(2,-1){\line(1,0){37}} \put(2,13){\line(1,0){37}}
%%These two are boundary lines

\put(6,5){\line(1,2){3}} \put(6,5){\line(-1,2){3}}

\put(9,6){\line(-6,5){1.5}} \put(3,11){\line(6,-5){1.1}}
\put(4.5,9.75){\line(6,-5){1.1}}
\put(6.15,8.375){\line(6,-5){0.88}} %(32.3,8.25)

\put(6.5,4.3){\mbox{$x_{ii}$}} \put(6.3,9.4){\mbox{$x_{ij}$}}
\put(3,5){\mbox{$x_{ji}$}} \put(3,11.5){\mbox{$x_{jj}$}}
\put(8,5){\mbox{$x_{ki}$}} \put(8,11.5){\mbox{$x_{kj}$}}

\put(3,3){\mbox{${x_ix_jx_k\times x_ix_j}$}}

%%%%%

%%%%%%%%%%%%%%%%%%%%%%%%

\thicklines

\put(13,5){\line(1,0){5}} \put(13,11){\line(1,0){5}}
\put(13,5){\line(1,3){1}} \put(13,5){\line(0,1){6}}
\put(14,8){\line(1,0){5}} \put(14,8){\line(-1,3){1}}
\put(18,5){\line(1,3){1}} \put(19,8){\line(-1,3){1}}
\put(18,5){\line(0,1){2.5}} \put(18,9){\line(0,1){2}}

\thinlines

\put(13,11){\line(2,-1){6}} \put(14,8){\line(4,-3){4}}
\put(13,11){\line(5,-6){1}} \put(14.25,9.5){\line(5,-6){1}}

\put(18,5){\line(-5,6){1}} \put(16.75,6.5){\line(-5,6){1}}

\put(13,11.4){\mbox{$x_{ii}$}} \put(17,11.4){\mbox{$x_{ji}$}}
\put(13,4){\mbox{$x_{ij}$}} \put(17,4){\mbox{$x_{jj}$}}
\put(13.9,6.7){\mbox{$x_{ik}$}} \put(19,7){\mbox{$x_{jk}$}}

\put(13,2.5){\mbox{${x_ix_j\times x_ix_jx_k}$}}

%%%%%

%%%%%%%%%%%%%%%%%%%%%%%%%%%%%%%%%%%%%%%%%%%%%%%%%%%%
\thicklines

\put(25,5){\line(-3,1){3}} \put(25,5){\line(3,1){3}}
\put(25,5){\line(0,1){5}} \put(22,6){\line(1,0){0.75}}
\put(23.25,6){\line(1,0){0.75}} \put(27.25,6){\line(1,0){0.75}}
\put(26,6){\line(1,0){0.75}} \put(22,6){\line(0,1){5}}
\put(28,6){\line(0,1){5}}

\put(25,10){\line(-3,1){3}} \put(25,10){\line(3,1){3}}
\put(22,11){\line(1,0){6}}

\thinlines

\put(25,5){\line(1,2){3}} \put(25,5){\line(-1,2){3}}

\put(28,6){\line(-6,5){1.5}} \put(22,11){\line(6,-5){1.1}}
\put(23.5,9.75){\line(6,-5){1.1}}
\put(25.15,8.375){\line(6,-5){0.88}} %(32.3,8.25)

\put(25.5,4.3){\mbox{$x_{ii}$}} \put(25.3,9.4){\mbox{$x_{il}$}}
\put(22,5){\mbox{$x_{ji}$}} \put(22,11.5){\mbox{$x_{j\,l}$}}
\put(27,5){\mbox{$x_{ki}$}} \put(27,11.5){\mbox{$x_{kl}$}}

\put(22,3){\mbox{${x_ix_jx_k\times x_ix_l}$}}

%%%%%

%%%%%%%%%%%%%%%%%%%%%%%%

\thicklines

\put(31,5){\line(1,0){5}} \put(31,11){\line(1,0){5}}
\put(31,5){\line(1,3){1}} \put(31,5){\line(0,1){6}}
\put(32,8){\line(1,0){5}} \put(32,8){\line(-1,3){1}}
\put(36,5){\line(1,3){1}} \put(37,8){\line(-1,3){1}}
\put(36,5){\line(0,1){2.5}} \put(36,9){\line(0,1){2}}

\thinlines

\put(31,11){\line(2,-1){6}} \put(32,8){\line(4,-3){4}}
\put(31,11){\line(5,-6){1}} \put(32.25,9.5){\line(5,-6){1}}

\put(36,5){\line(-5,6){1}} \put(34.75,6.5){\line(-5,6){1}}

\put(31,11.4){\mbox{$x_{ii}$}} \put(35,11.4){\mbox{$x_{li}$}}
\put(31,4){\mbox{$x_{ij}$}} \put(35,4){\mbox{$x_{lj}$}}
\put(31.9,6.7){\mbox{$x_{ik}$}} \put(37,7){\mbox{$x_{lk}$}}

\put(31,2.5){\mbox{${x_ix_l\times x_ix_jx_k}$}}

%%%%%

%%%%%%

\put(5,1.2){\boldmath $(a)$} \put(15,1.2){\boldmath $(b)$}
\put(24.3,1.2){\boldmath $(c)$} \put(33,1.2){\boldmath $(d)$}

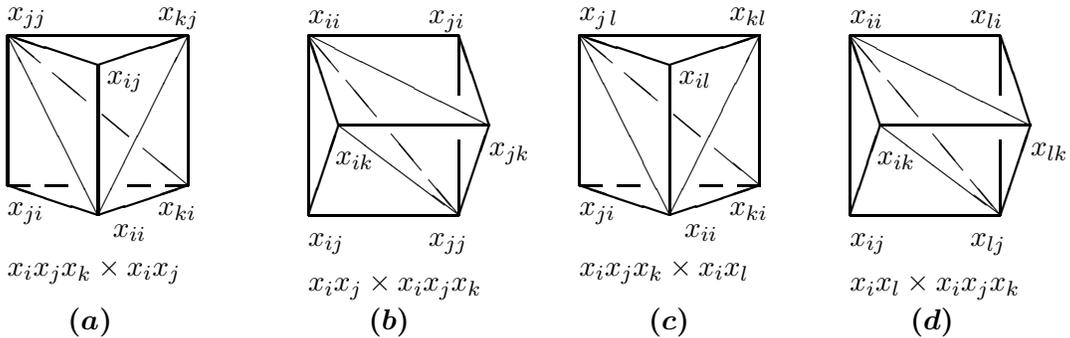
\put(6.5,-0.3){\mbox{\bf Figure 6\,: Simplicial subdivisions of
3-cells of {\boldmath $S^{\,2}_4 \times S^{\,2}_4$}}}

\end{picture}

\medskip

For $i, j, k, l$ distinct, consider the 3-cell $x_ix_jx_k \times
x_ix_l$. Here $x_{ii}x_{j\,l}$ and $x_{ii}x_{kl}$ are edges. By
interchanging $j$ and $k$ (if required) we may assume that $\{i,
j, k, l\}$ is an even permutation of $\{1, 2, 3, 4\}$. Then
$x_{ki}x_{j\,l}$ is an edge and hence, by Lemma \ref{L2},
$x_ix_jx_k \times x_ix_l = x_{ii}x_{il}x_{kl}x_{j\,l} \cup
x_{ii}x_{ki}x_{kl}x_{j\,l} \cup x_{ii}x_{ki}x_{ji}x_{j\,l}$ is the
unique subdivision of $x_ix_jx_k \times x_ix_l$ (cf. Figure 6\,$(c)$). Similarly, for
the 3-cell $x_ix_l\times x_ix_jx_k$, we may assume that $\{i, j,
k, l\}$ is an even permutation of $\{1, 2, 3, 4\}$. Then
$x_{ik}x_{lj}$ is an edge and hence, by Lemma \ref{L2}, $x_ix_l
\times x_ix_jx_k  = x_{ii}x_{li}x_{lk}x_{lj} \cup x_{ii}x_{ik}
x_{lk}x_{lj} \cup x_{ii}x_{ik}x_{ij}x_{lj}$ is the unique
subdivision of $x_ix_l \times x_ix_jx_k$ (cf. Figure 6\,$(d)$). These give the
3-skeleton of $X$.

For $i, j, k$ distinct, consider the 4-cell $A= x_ix_jx_k \times
x_ix_jx_k$. The boundary $\partial A$ of $A$ consists of six
3-cells. From above, it follows that $S^1(\{x_{ii}, x_{jj},
x_{kk}\}) \ast C_6 \subseteq X$ is the subdivision of $\partial
A$, where $C_6$ is the 6-cycle $C_6(x_{ij}, x_{ik}, x_{jk},
x_{ji}, x_{ki}, x_{kj})$. Let $D \subseteq X$ be the subdivision
of $A$. Then, $D$ is a 9-vertex 4-ball with boundary $\partial D =
S^1(\{x_{ii}, x_{jj}, x_{kk}\}) \ast C_6$. Clearly, $C_6$ is an
induced subcomplex of $X$. Therefore, each 4-simplex in $B$ must
contain $x_{ii}x_{jj}x_{kk}$. Thus, $x_{ii}x_{jj}x_{kk}$ is a
simplex in $D \setminus \partial D$. Therefore, ${\rm
lk}_D(x_{ii}x_{jj}x_{kk})$ is a cycle and hence $= C_6$. These
imply that $D = \overline{x_{ii}x_{jj}x_{kk}}\ast C_6$.

Now, consider the 4-cell $B = x_ix_jx_k \times x_ix_jx_l$, where
$i, j, k, l$ are distinct. By interchanging $i$ and $j$ (if
required) we may assume that $\{i, j, k, l\}$ is an even
permutation of $\{1, 2, 3, 4\}$. The boundary $\partial B$ of $B$
consists of six 3-cells. From above, it follows that the
subdivision of $\partial B$ in $X$ is a 9-vertex triangulated
3-sphere and obtained from $S^1_3(\{x_{ii}, x_{jj}, x_{kl}\})
\times C_5$ by starring the vertex $x_{ji}$ in the 3-simplex
$\alpha := x_{ii} x_{jj}x_{jl}x_{ki}$, where $C_5$ is the 5-cycle
$C_5(x_{ij}, x_{il}, x_{jl}, x_{ki}, x_{kj})$. Since $x_{ji}
x_{ij}$, $x_{ji} x_{il}$, $x_{ji}x_{kj}$ and $x_{ji} x_{kl}$ are
non-edges, it follows that $\sigma := x_{ii}x_{jj} x_{ji} x_{jl}
x_{ki}$ is the only possible 4-simplex containing $x_{ji}$ inside
$B$. So, $\sigma \in X$. Then $B = \sigma \cup P$, where $P$ is a
4-cell such that $P \cap \sigma = \alpha$ and $S^1_3(\{x_{ii},
x_{jj}, x_{kl}\}) \ast C_5 \subseteq X$ is the subdivision of
$\partial P$ in $X$ (i.e., $P$ is the 4-cell whose geometric
carrier is $(|B|\setminus |\sigma|)\cup |\alpha|$). Let $Q$ be the
simplicial subdivision of $P$ in $X$. So, $\partial Q=
S^1_3(\{x_{ii}, x_{jj}, x_{kl}\}) \ast C_5$. Since $C_5$ is
induced in $X$, it follows that any 4-simplex in $Q$ must contain
$x_{ii}x_{jj} x_{kl}$. Since $x_{ii}x_{jj} x_{kl}\in Q \setminus
\partial Q$, ${\rm lk}_Q(x_{ii}x_{jj} x_{kl})$ is a cycle and
hence $= C_5$. These imply that $Q = \overline{x_{ii} x_{jj}
x_{kk}} \ast C_5$. Then $B = (\overline{x_{ii}x_{jj} x_{kl}} \ast
C_5) \cup \bar{\sigma}$.

Now, we have subdivided all the 4-cells in $S^{2}_4 \times
S^{2}_4$. It is routine to check that the resulting simplicial
complex $X$ is identical with the complex $(S^{2}\times
S^{2})_{16}$ defined in Section 1.

\medskip

\noindent {\bf Case 2.} $x_1x_3 \times x_2x_4 = x_{12}x_{32}x_{14}
\cup x_{32}x_{14}x_{34}$. By the same method as in Case 1, one can
show that $X$ is uniquely determined and is isomorphic to
$(S^{\,2} \times S^{\,2})_{16}$ via the map $f$ given by the
transposition $(1,2)$ on the suffixes, i.e., $f \equiv (x_{11}
x_{22})(x_{13} x_{23}) (x_{14} x_{24})(x_{31} x_{32})(x_{41}
x_{42})$. This completes the proof. \hfill $\Box$

\bigskip

\noindent {\bf Proof of Corollary \ref{C2}.}
From Proposition \ref{P6}, Lemma \ref{L1} and  Theorem \ref{T1},
it is immediate that $\CC P^{\,2}_{10}$ triangulates $\CC P^{\,2}$.

Since the automorphism group $A_4 = \langle\alpha, \beta\rangle$
of $(S^{\,2}\times S^{\,2})_{16}$ commutes with $\ZZ_2$, it
descends  to an automorphism group $A_4 = \langle \bar{\alpha},
\bar{\beta} \rangle$ of $\CC P^{\,2}_{10}$. We need to show that
there are no other automorphisms.

It is easy to check that the four vertices $x_{ii}$, $1\leq i\leq
4$, are the only ones with 2-neighborly links. Therefore, the full
automorphism group must fix this set of four vertices. Since $A_4$
is 2-transitive on this 4-set, it suffices to show that there is
no non-trivial automorphism $\gamma$ fixing both $x_{11}$ and
$x_{22}$. Suppose the contrary. Then $\gamma$ is a non-trivial
automorphism of ${\rm lk}(x_{11}x_{22})$. But ${\rm
lk}(x_{11}x_{22})$ is the 8-vertex triangulated 2-sphere given in
Figure 7.

\bigskip

\setlength{\unitlength}{3mm}

\begin{picture}(48,16.5)(-3,-3.5)

%%%%%%%%%%%%%%%%%%%%%%% %%%%%%%%%%%%%%%%%%%%%%%%%%%%%
\thicklines

\put(0,0){\line(1,0){24}} \put(0,0){\line(1,1){12}}
\put(0,0){\line(6,1){12}} \put(0,0){\line(3,1){9}}
\put(0,0){\line(5,3){10}}

\put(24,0){\line(-1,1){12}} \put(24,0){\line(-6,1){12}}
\put(24,0){\line(-3,1){9}} \put(24,0){\line(-5,3){10}}

\put(12,2){\line(3,1){3}} \put(12,2){\line(1,2){2}}
\put(12,2){\line(-3,1){3}} \put(12,2){\line(-1,2){2}}
\put(12,2){\line(0,1){10}}

\put(9,3){\line(1,3){3}} \put(15,3){\line(-1,3){3}}

\put(-1.5,1){\mbox{$x_{24}$}} \put(23.8,1){\mbox{$x_{13}$}}
\put(11.2,0.8){\mbox{$x_{12}$}} \put(13,11.6){\mbox{$x_{34}$}}

\put(7.2,3.4){\mbox{$x_{44}$}} \put(15.3,3.4){\mbox{$x_{33}$}}
\put(8.2,6.4){\mbox{$x_{14}$}} \put(14.3,6.4){\mbox{$x_{23}$}}

\put(5,-2.5){\mbox{\bf Figure 7\,: ${\rm lk}_{\C
P^2_{10}}(x_{11}x_{22})$}}

\thinlines

\put(-4,-3.5){\line(1,0){49}} \put(-4,13){\line(1,0){49}}
%%These two are boundary lines

\end{picture}

\medskip

From the picture, it is apparent that ${\rm lk}(x_{11}x_{22})$ has
only one non-trivial automorphism, namely $(x_{13}, x_{24})
(x_{14}, x_{23}) (x_{33}, x_{44})$. Therefore, $\gamma = (x_{13},
x_{24}) (x_{14}, x_{23}) (x_{33}, x_{44})$ and hence $\gamma$
fixes the 3-simplex $x_{11}x_{33}x_{44}x_{34}$. Then $\gamma$ must
either fix or interchange the two vertices $x_{13}$ and $x_{14}$
in the link of this 3-simplex, a contradiction. This completes the
proof. \hfill $\Box$

\bigskip

\noindent {\bf Proof of Theorem \ref{T3}.} Consider the following
sequence of bistellar moves on $\CC P^{\,2}_{10}$ (performed
one after the other)\,:
\begin{eqnarray*}
\mbox{(i)} ~ x_{22}x_{33}x_{44} \mapsto x_{23}x_{24}x_{34}, &
\mbox{(ii)} ~ x_{11}x_{33}x_{44} \mapsto x_{13}x_{14}x_{34}, &
\mbox{(iii)} ~ x_{11}x_{22}x_{44} \mapsto x_{12}x_{14}x_{24}, \\
\mbox{(iv)} ~  x_{14}x_{33}x_{44} \mapsto x_{12}x_{13}x_{34}, &
\mbox{(v)} ~ x_{22}x_{34}x_{44} \mapsto x_{13}x_{23}x_{24}, &
\mbox{(vi)} ~ x_{23}x_{33}x_{44} \mapsto x_{12}x_{24}x_{34}, \\
\mbox{(vii)} ~ x_{12}x_{22}x_{44} \mapsto x_{13}x_{14}x_{24}, &
\mbox{(viii)} ~ x_{33}x_{44} \mapsto x_{12}x_{13}x_{24}x_{34}, &
\mbox{(ix)} ~ x_{22}x_{44} \mapsto x_{13}x_{14}x_{23}x_{24}.
\end{eqnarray*}
At the end of these moves, we get a 10-vertex triangulation
$K$ of $\CC P^{\,2}$.

On $K$ we perform the following sequence of bistellar moves one
after another.
\begin{eqnarray*}
\mbox{(x)} ~ x_{11}x_{24}x_{44} \mapsto x_{12}x_{14}x_{23}, \!\!&\!\!
\mbox{(xi)} ~ x_{11}x_{13}x_{44} \mapsto x_{14}x_{23}x_{34}, \!\!&\!\!
\mbox{(xii)} ~ x_{11}x_{44} \mapsto x_{12}x_{14}x_{23}x_{34}, \\
\mbox{(xiii)} ~ x_{44}x_{14}x_{24} \mapsto x_{12}x_{13}x_{23}, \!\!&\!\!
\mbox{(xiv)} ~ x_{44}x_{14} \mapsto x_{34}x_{12}x_{13}x_{23}, \!\!&\!\!
\mbox{(xv)} ~ x_{44} \mapsto x_{24}x_{34}x_{12}x_{13}x_{23}.
\end{eqnarray*}
(Note that the last three bistellar moves together is same as
the GBM with respect to $(\overline{x_{44}} \ast S^1_3(\{x_{14},
x_{24}, x_{34}\}) \ast S^1_3(\{x_{12}, x_{13}, x_{23}\}),
S^1_3(\{x_{14}, x_{24}, x_{34}\}) \ast \overline{x_{12}x_{13}x_{23}})$.)
The last bistellar move deletes the vertex $x_{44}$
and hence obtain a 9-vertex triangulation $L$ of $\CC P^{\,2}$.
(Observe that $A_1= \{x_{11}, x_{23}, x_{24}\}$, $A_2= \{x_{14},
x_{33}, x_{12}\}$, $A_3= \{x_{34}, x_{22}, x_{13}\}$ is an
amicable partition of $L$ whose layer is of first type (cf.
\cite{bd3}).)

Let $\CC P^{\,2}_9$ be as described in \cite{kb2} with vertex-set
$\{1, 2, \dots, 9\}$. Consider the map $\varphi \colon L \to \CC
P^{2}_9$ given by: $\varphi(x_{11}) = 1$, $\varphi(x_{23}) = 2$,
$\varphi(x_{24}) = 3$, $\varphi(x_{34}) = 4$, $\varphi(x_{22}) =
5$, $\varphi(x_{13}) = 6$, $\varphi(x_{14}) = 7$, $\varphi(x_{33})
= 8$, $\varphi(x_{12}) = 9$. It is easy to see that $\varphi$ is
an isomorphism. Thus, $\CC P^{\,2}_{10}$ is bistellar equivalent to
$\CC P^{\,2}_9$.

\medskip

Now, on $K$ we perform the following sequence of bistellar moves\,:
\begin{eqnarray*}
\mbox{(xvi)} ~ x_{11}x_{22}x_{33} \mapsto x_{12}x_{13}x_{23}, &
\mbox{(xvii)} ~ x_{22}x_{33}x_{24} \mapsto x_{14}x_{23}x_{34}, \\
\mbox{(xviii)} ~ x_{22}x_{33}x_{13} \mapsto x_{12}x_{14}x_{23}, &
\mbox{(xix)} ~ x_{22}x_{33} \mapsto x_{12}x_{14}x_{23}x_{34}.
\end{eqnarray*}
We obtain a 10-vertex triangulation $M$ of $\CC P^{\,2}$.
Let $K^4_{10}$ be as described in \cite{kb1} with vertex-set
$\{X, Y, Z, 0, 1, \dots, 6\}$. Consider the map $\psi \colon
M \to K^4_{10}$ given by $\psi(x_{33}) = X$, $\psi(x_{22}) = Y$,
$\psi(x_{44}) = Z$, $\psi(x_{11}) = 0$, $\psi(x_{13}) = 1$,
$\psi(x_{12}) = 2$, $\psi(x_{23}) = 3$, $\psi(x_{14}) = 4$,
$\psi(x_{34}) = 5$, $\psi(x_{24}) = 6$. It is easy to see that
$\psi$ is an isomorphism. Thus, $\CC P^{\,2}_{10}$ is bistellar
equivalent to $K^{4}_{10}$. This completes the proof. \hfill $\Box$

\bigskip

\noindent {\bf Proof of Corollary \ref{C4}.} Part $(a)$ follows from
Corollary \ref{C2} and Theorem \ref{T3}.

In \cite{kb1}, explicit coordinates for simplices of $K^4_{10}$ in
the Fubini-Study metric were given. This shows that the induced
pl-structure on $\CC P^2$ by $K^4_{10}$ is the standard one. Part
$(b)$ now follows from Theorem \ref{T3}. \hfill $\Box$

\begin{lemma}$\!\!\!${\bf .} \label{L3}
Let $I_1$ and $I_2$ be an antimorphic pair of icosahedra. Then we
have\,:  \vspace{-2mm}
\begin{enumerate}
\item[{\rm (a)}] ${\rm Aut}(I_1) = {\rm Aut}(I_2) = A_5\times
\ZZ_2$. \vspace{-2mm}
 \item[{\rm (b)}] For each triangle $\Delta$ of $I_1$, there is
a unique triangle $\Delta^{\prime}$ of $I_2$ such that each of the
three triangles of $I_2$ sharing an edge with $\Delta^{\prime}$
has its third vertex in $\Delta$. Further, the map $\varphi \colon
\Delta \mapsto \Delta^{\prime}$ is a bijection from the triangles
of $I_1$ to the triangles of $I_2$. There is a similarly defined
bijection $\psi$ from the triangles of $I_2$ to the triangles of
$I_1$, and \vspace{-2mm}
 \item[{\rm (c)}] Every isomorphism $f \colon I_1 \to I_2$
 intertwines $\,\varphi$ and $\psi$.
\end{enumerate}
\end{lemma}

(Warning\,: The maps $\varphi$ and $\psi$ are {\bf not} induced by
any vertex - to - vertex map\,!)

\bigskip

\noindent {\bf Proof.} Recall that $I_1$ and $I_2$ have the same
vertex set and the same pairs of antipodal vertices. Thus, they
have the same antipodal map (sending each vertex $x$ to its
antipode $\bar{x}$). Now, the full automorphism group of the
icosahedron is generated by its rotation group $A_5$ and the
antipodal map. So, to prove Part (a), it suffices to show that
$I_1$ and $I_2$ share the same rotation group. For each pair $x$,
$\bar{x}$ of antipodes, $I_i$ has a rotation symmetry
$\alpha^i_{x, \bar{x}}$ which fixes $x$ and $\bar{x}$ and rotates
the remaining vertices along the 5-cycles ${\rm lk}_{I_i}(x)$ and
${\rm lk}_{I_i}(\bar{x})$. The rotation group of $I_i$ is
generated by these automorphisms of order five. But, ${\rm
lk}_{I_2}(x)$ (respectively, ${\rm lk}_{I_2}(\bar{x})$) is the
graph theoretic complement of the pentagon ${\rm
lk}_{I_1}(\bar{x})$ (respectively, ${\rm lk}_{I_1}(x)$).
Therefore, $\alpha^2_{x, \bar{x}}$ is the square of $\alpha^1_{x,
\bar{x}}$. This proves Part (a).

Notice that if $f_1, f_2 \colon I_1 \to I_2$ are two
antimorphisms, then $f_1 \circ f_2^{-1} \in {\rm Aut}(I_2)$ and
$f_2^{-1} \circ f_1 \in {\rm Aut}(I_1)$. Thus, the antimorphism is
unique up to right multiplication by elements of ${\rm Aut}(I_1)$
(or left multiplication by elements of ${\rm Aut}(I_2)$).
Therefore, there is no loss of generality in taking the
antimorphic pair of icosahedra as the one given in Figure 2.

Since the common automorphism group is transitive on the triangles
of $I_1$ (and of $I_2$), it is enough to look at the triangle
$\Delta = x_{12}x_{13}x_{14}$ of $I_1$. From Figure 2, we see that
the links in $I_2$ of two vertices of $\Delta$ have exactly two
vertices in common. Namely, we have $V({\rm lk}_{I_2}(x_{12}))
\cap V({\rm lk}_{I_2}(x_{13})) = \{x_{21}, x_{32}\}$, $V({\rm
lk}_{I_2}(x_{12})) \cap V({\rm lk}_{I_2}(x_{14})) = \{x_{24},
x_{41}\}$, $V({\rm lk}_{I_2}(x_{13})) \cap V({\rm
lk}_{I_2}(x_{14})) = \{x_{31}, x_{43}\}$. Therefore, any triangle
$\Delta^{\prime}$ of $I_2$ satisfying the requirement must be
contained in the vertex set $\{x_{21}, x_{32}, x_{24}, x_{41},
x_{31}, x_{43}\}$. But one sees that this set of six vertices
contains a unique triangle in $I_2$, namely $\Delta^{\prime} =
x_{21}x_{31}x_{41}$. Thus the map $\varphi \colon \Delta \to
\Delta^{\prime}$ is well defined. Similarly, there is a well
defined map $\psi$ from the triangles of $I_2$ to the triangles of
$I_1$. The map $\psi \circ \varphi$ is the antipodal map on the
triangles of $I_1$ to themselves. Similarly, $\varphi \circ \psi$
is the antipodal map on triangles of $I_2$. Hence $\varphi$ (as
well as $\psi$) is a bijection. This proves Part (b).

To prove Part (c), let $f$ be any isomorphism from $I_1$ to $I_2$.
Since $I_1$ and $I_2$ are antimorphic, it is immediate that $f$
also defines an isomorphism from $I_2$ to $I_1$. Let $\Delta$ be
any triangle of $I_1$ and let $\Delta^{\prime} = \varphi(\Delta)$.
By definition, there are three triangles $\Delta^{\prime}_1$,
$\Delta^{\prime}_2$, $\Delta^{\prime}_3$ of $I_2$ each of which
shares a vertex with $\Delta$ and an edge with $\Delta^{\prime}$.
Then $f(\Delta)$ and $f(\Delta^{\prime})$ are triangles of $I_2$
and $I_1$, respectively. Also, $f(\Delta^{\prime}_1)$,
$f(\Delta^{\prime}_2)$, $f(\Delta^{\prime}_3)$ are three triangles
of $I_1$ each of which shares a vertex with $f(\Delta)$ and an
edge with $f(\Delta^{\prime})$. Therefore, by definition of
$\psi$, $\psi(f(\Delta))= f(\Delta^{\prime})= f(\varphi(\Delta))$.
\hfill $\Box$

\bigskip

\noindent {\bf Proof of Theorem \ref{T5}.} As in the proof of
Theorem \ref{T1}, one
may verify that $(S^{\,2} \times S^{\,2})_{16}^{\prime}$ is a
simplicial subdivision of $S^{\,2}_4 \times S^{\,2}_4$, and hence
it triangulates $S^{\,2} \times S^{\,2}$. We apply the following
sequence of bistellar moves to $(S^{\,2} \times
S^{\,2})_{16}^{\prime}$ to create a second 16-vertex
triangulation $(S^{\,2}\times \widehat{S^{\,2})_{16}}$ of $S^{\,2}
\times S^{\,2}$\,: \vspace{-1mm}
\begin{eqnarray*}
x_{12}x_{13}x_{14} \mapsto x_{23}x_{34}x_{42}, &
x_{21}x_{23}x_{24} \mapsto x_{14}x_{31}x_{43}, \\
x_{31}x_{32}x_{34} \mapsto x_{12}x_{24}x_{41}, &
x_{41}x_{42}x_{43} \mapsto x_{13}x_{21}x_{32}.
\end{eqnarray*}
Since this set of bistellar moves is stable under the automorphism
group $A_4$ of $(S^{\,2} \times S^{\,2})_{16}^{\prime}$, it
follows that $(S^{\,2} \times \widehat{S^{\,2})_{16}}$ inherits
the group $A_4$. Also, both complexes have ${\rm lk}(x_{11}) =
S^1_3(\{x_{12}, x_{13}, x_{14}\}) \ast S^1_3(\{x_{21}, x_{31},
x_{41}\})$. However, while $(S^{\,2} \times
S^{\,2})_{16}^{\prime}$ has both $x_{12}x_{13}x_{14}$ and
$x_{21}x_{31}x_{41}$ as triangles, we have chosen the bistellar
moves judiciously to ensure that $(S^{\,2} \times
\widehat{S^{\,2})_{16}}$ does not have the triangle $x_{12}x_{13}
x_{14}$. Therefore, we may apply the following four GBM's (one
after the other) to $(S^{\,2} \times \widehat{S^{\,2})_{16}}$ to
delete the four vertices $x_{ii}$, $1\leq i\leq 4$\,:
\begin{eqnarray*}
({\rm st}(x_{11}), D^2_3(\{x_{12}, x_{13}, x_{14}\}) \ast
S^1_3(\{x_{21}, x_{31}, x_{41}\}), \\
({\rm st}(x_{22}), D^2_3(\{x_{21}, x_{23}, x_{24}\}) \ast
S^1_3(\{x_{12}, x_{32}, x_{42}\}), \\
({\rm st}(x_{33}), D^2_3(\{x_{31}, x_{32}, x_{34}\}) \ast
S^1_3(\{x_{13}, x_{23}, x_{43}\}), \\
({\rm st}(x_{44}), D^2_3(\{x_{41}, x_{42}, x_{43}\}) \ast
S^1_3(\{x_{14}, x_{24}, x_{34}\}).
\end{eqnarray*}
The resulting complex $X$ is therefore a 12-vertex triangulation
of $S^{\,2} \times S^{\,2}$. So, to confirm the first statement of
this theorem, it suffices to show that $X$ is isomorphic to the
complex $(S^{\,2}\times S^{\,2})_{12}$ described in Section 3.
Indeed, with the antimorphic pair of icosahedra (and their vertex
names) as in Figure 2, we shall show that we actually have $X =
(S^{\,2}\times S^{\,2})_{12}$.

Notice that $X$ inherits the automorphism group $A_4$ from
$(S^{\,2} \times \widehat{S^{\,2})_{16}}$, and modulo this group,
the following six are basic facets of $X$\,: \vspace{-1mm}
\begin{eqnarray*}
x_{12}x_{14}x_{21}x_{24}x_{31}, & x_{12}x_{13}x_{14}x_{21}x_{31},
& x_{12}x_{23}x_{31}x_{13}x_{32}, \\
x_{12}x_{31}x_{34}x_{14}x_{24}, & x_{24}x_{31}x_{32}x_{12}x_{21},
& x_{24}x_{31}x_{32}x_{12}x_{41}.
\end{eqnarray*}
Each basic facet is in an $A_4$-orbit of size 12, yielding a total
of $6 \times 12 = 72$ facets of $X$. Since $(S^{\,2} \times
S^{\,2})_{12}$ also has 72 facets and since the group $A_4$
(acting on subscripts) is a subgroup of the automorphism group
$A_5 \times \ZZ_2$ of $(S^{\,2}\times S^{\,2})_{12}$, it suffices
to observe that all six basic facets of $X$ listed above are also
facets of $(S^{\,2}\times S^{\,2})_{12}$. Indeed, the first facet
$x_{12} x_{14}x_{21}x_{24}x_{31}$ is in $(S^{\,2}\times
S^{\,2})_{12}$ since these five vertices are the neighbors of
$x_{23}$ in $I_1$ (and of $x_{41}$ in $I_2$). In each of the
remaining five basic facets of $X$, the first three vertices
constitute a triangle $\Delta$ of $I_1$ with the last two vertices
in the corresponding triangle $\Delta^{\prime} = \varphi(\Delta)$
of $I_2$ (cf. Lemma \ref{L3}). (For instance, $\Delta =
x_{12}x_{13}x_{14}$ is a triangle of $I_1$, with corresponding
triangle $\Delta^{\prime} = x_{21}x_{31}x_{41}$ of $I_2$.
Therefore, the second basic facet of $X$ is a facet of
$(S^{\,2}\times S^{\,2})_{12}$.) This shows that $(S^{\,2}\times
S^{\,2})_{12} = X$, so that $(S^{\,2}\times S^{\,2})_{12}$
triangulates $S^{\,2}\times S^{\,2}$.

To compute the full automorphism group of $(S^{\,2}\times
S^{\,2})_{12}$, notice that it has exactly 40 triangles of degree
3 (the rest are of degree 5), namely the twenty triangles of $I_1$
and the twenty triangles of $I_2$. Consider the graph whose
vertices are these forty triangles, two of them being adjacent if
and only if they share an edge. This graph has exactly two
connected components, of size 20 each, namely the triangles of
$I_1$ and $I_2$. This shows that any automorphism $f$ of
$(S^{\,2}\times S^{\,2})_{12}$ either fixes both $I_1$ and $I_2$
or interchanges them. So, ${\rm Aut}(I_1) = {\rm Aut}(I_2)= A_5
\times \ZZ_2$ is a subgroup of index at most two in the full
automorphism group of $(S^{\,2}\times S^{\,2})_{12}$.

Let $f\colon I_1 \to I_2$ be any isomorphism. Since $I_1$ and
$I_2$ are antimorphic, it is immediate that $f$ is also an
isomorphism from $I_2$ to $I_1$. Since the five neighbors in $I_1$
of any vertex are also the neighbors in $I_2$ of the antipodal
vertex, it is immediate that $f$ maps each of the 12 facets of the
first kind in $(S^{\,2}\times S^{\,2})_{12}$ to a facet of the
same kind. Also, for any triangle $\Delta$ of $I_1$, the
construction of $(S^{\,2}\times S^{\,2})_{12}$ shows that ${\rm
lk}(\Delta) = S^1_3(\varphi(\Delta))$, and also, for any triangle
$\Delta^{\prime}$ of $I_2$, ${\rm lk}(\Delta^{\prime}) =
S^1_3(\psi(\Delta^{\prime}))$. Since $f$ intertwines $\varphi$ and
$\psi$ (Lemma \ref{L3}), we also have ${\rm lk}(f(\Delta)) =
S^1_3(\psi(f(\Delta))) = S^1_3(f(\varphi(\Delta))) =
f(S^1_3(\varphi(\Delta))) = f({\rm lk}(\Delta))$. Similarly, for
any triangle $\Delta^{\prime}$ of $I_2$, ${\rm
lk}(f(\Delta^{\prime})) = f({\rm lk}(\Delta^{\prime}))$. Thus, $f$
also maps all sixty facets of the second type in $(S^{\,2}\times
S^{\,2})_{12}$ to facets of the same type. Thus, any isomorphism
between $I_1$ and $I_2$ is also an automorphism of $(S^{\,2}\times
S^{\,2})_{12}$. Therefore, the full automorphism group $G$ of
$(S^{\,2}\times S^{\,2})_{12}$ has $H = A_5 \times \ZZ_2$ as an
index two subgroup. Thus, $G$ is of order 240. Indeed, $G$
consists of the 120 common automorphisms of $I_1$ and $I_2$, and
the 120 isomorphisms between $I_1$ and $I_2$. In particular, take
$g = (x_{12} x_{21} x_{24} x_{42} x_{14} x_{41} x_{43} x_{34}
x_{13} x_{31} x_{32} x_{23})$, which is an isomorphism between
$I_1$ and $I_2$. Note that $g^6$ is the common antipodal map of
$I_1$ and $I_2$, hence it is in the center of $G$. Thus,
$G/\langle g^6\rangle$ is the extension of $A_5$ by the involution
$\alpha = g$ (mod $g^6$). But $A_5$ has only one non-trivial
extension by an involution, namely $S_5$. So, $G$ is an extension
of a central involution by $S_5$. It can not be the split
extension $S_5 \times \ZZ_2$ since this has no element of order
12. Therefore, $G$ is the unique non-split extension $2S_5$ of
$\ZZ_2$ by $S_5$. \hfill $\Box$

\begin{remark}$\!\!\!${\bf .} \label{R5}
{\rm If the link of a vertex $u$ in a triangulated 4-manifold $X$
is $S^1_3(\{x, y, z\})\ast S^1_3(\{a, b, c\})$ and $xyz$ is not a
simplex in $X$ then the GBM $({\rm st}_X(u), D^2_3(\{x, y, z\})
\ast S^1_3(\{a, b, c\})$ is equivalent to the sequence of the
following three bistellar moves\,: $uab \mapsto xyz$, $ua \mapsto
cxyz$, $u \mapsto bcxyz$. Thus, from the proof of Theorem
\ref{T5}, $(S^{\,2}\times S^{\,2})_{12}$ can be obtained from
$(S^{\,2}\times S^{\,2})_{16}^{\prime}$ by a sequence of bistellar
moves only. }
\end{remark}

\medskip

\noindent {\bf Acknowledgement.} The authors thank the anonymous
referee for many useful comments which led to substantial
improvements in the presentation of this paper. The authors are
thankful to Siddhartha Gadgil, Frank Lutz, Wolfgang K\"{u}hnel and
Alberto Verjovsky for useful conversations and references. The
authors thank Ipsita Datta for her help in the proof of Theorem
\ref{T3}.
%The second author was partially supported by UGC-SAP/DSA-IV.

%\newpage

\end{document}